\RequirePackage{etex}
\documentclass{article}
\usepackage{graphicx} 
\usepackage{amsmath}
\usepackage{subfig}
\graphicspath{ {./img/} }
\usepackage{autonum}
\usepackage{float}
\usepackage{svg}
 
\newtheorem{remark}{Remark}

\makeatletter
\let\@fnsymbol\@arabic
\makeatother

\title{Unconditionally local bounds preserving numerical scheme based on inverse Lax-Wendroff procedure for advection on networks}
\author{Peter Frolkovi\v{c} \thanks{peter.frolkovic@stuba.sk; Supported by the VEGA grant 1/0314/23 and APVV-23-0186.}\and Svetlana Kri\v{s}ková \thanks{svetlana.kriskova@stuba.sk; Supported by the VEGA grant 1/0314/23.} \and Katar\'ina Lackov\'a \thanks{katarina.lackova@stuba.sk; Funded by the EU NextGenerationEU through the Recovery and Resilience Plan for Slovakia under the project No. 09I03-03-V05-00005.}}

\date{\today}

\begin{document}

\maketitle

\begin{abstract}
We derive an implicit numerical scheme for the solution of advection equation where the roles of space and time variables are exchanged using the inverse Lax-Wendroff procedure. 
The scheme contains a linear weight for which it is always second order accurate in time and space, and the stencil in the implicit part is fully upwinded for any value of the weight, enabling a direct computation of numerical solutions by forward substitution.  
To fulfill the local bounds for the solution represented by the discrete minimum and maximum principle (DMP), we use a predicted value obtained with the linear weight and check a priori if the DMP is valid. If not, we can use either a nonlinear weight or a limiter function that depends on Courant number and apply such a high-resolution version of the scheme to obtain a corrected value. 
The advantage of the scheme obtained with the inverse Lax-Wendroff procedure is that only in the case of too small Courant numbers, the limiting is towards the first order accurate scheme, which is not a situation occurring in numerical simulations with implicit schemes very often. In summary, the local bounds are satisfied up to rounding errors unconditionally for any Courant numbers, and the formulas for the predictor and the corrector are explicit. The high-resolution scheme can be extended straightforwardly for advection with nonlinear retardation coefficient with numerical solutions satisfying the DMP, and a scalar nonlinear algebraic equation has to be solved to obtain each predicted and corrected value. In numerical experiments, including transport on a sewer network, we can confirm the advantageous properties of numerical solutions for several representative examples.
\end{abstract}

\section{Introduction}

Numerical simulations can be a useful tool in a study of many practical applications; therefore, in engineering practice, there is a demand for fast and reliable software libraries for such a purpose. If a mathematical model behind the studied application is given by partial differential equations, numerical methods must be implemented if one does not want to restrict only to simplified applications. 
Potential users of such numerical tools have clear demands on their efficiency and reliability, and they expect no restrictions based solely on requirements of numerical methods that are not grounded by the application. However, such conditions are often present in numerical tools in the form of constrains on the choice of parameters in numerical methods. Clearly, users must choose appropriate parameters to reach satisfactory accuracy of approximate solutions, but, unfortunately, such choices need not to ensure physically acceptable numerical solutions or even stable computations with no overflow errors in computations.

When taking into account such demands of users in the case of numerical simulations for transport problems on networks, one can recognize these difficulties already for transient simulations of the advection equation, which is very often the dominant feature of more involved mathematical models. We are interested here mainly in transport problems in sewage systems \cite{markRiskAnalysesSewer1998,banik2014swmm5,veliskovaInverseTaskPollution2023}, but the advection-dominated problems also arise in other contexts, such as district heating networks \cite{eimerImplicitFiniteVolume2022,mohring2021district}. 

In deciding between two basic classes of time discretization to work with, we choose implicit time discretization methods. The main motivation to use it is to completely remove any stability restrictions on numerical discretization steps, so the choice of time steps can be dictated purely by accuracy requirements. Such a property is attractive for all users of simulation software in general, but for many cases like a large variation of characteristic speeds, it might even be necessary for practical reasons. In the case of networks, it means, for instance, that a single short segment (an edge) of the network does not dictate the choice of very small discretization steps to preserve stability of computations, whereas accuracy requirements could allow much larger steps and, consequently, much faster computations.

The main disadvantage of implicit discretizations is the necessity of solving algebraic equations to obtain the numerical solution, whereas explicit numerical schemes define its values directly. This disadvantage does not hold for algebraic systems with a simpler structure that can be obtained for special types of flow and transport problems, such as one-dimensional problems for which the characteristic speed does not change sign. This property is used, for example, in \cite{eimerImplicitFiniteVolume2022,barsukowImplicitActiveFlux2024} for advection in networks with a predetermined direction of flow, and in our study we also rely on this property. In \cite{eimerImplicitFiniteVolume2022,barsukowImplicitActiveFlux2024}, implicit numerical schemes from the third order up to the fifth order of accuracy are successfully proposed that are constructed purely on the upwind principle, producing systems of algebraic equations with a lower triangular matrix. Consequently, the system can be solved directly by a forward substitution that makes these implicit methods competitive with respect to the explicit ones, especially if the stability of the implicit method is taken into account.  

Unfortunately, the availability of efficient algebraic solvers is not enough in general when using implicit methods to solve flow and transport problems for real applications. It is known that when applying higher-order accurate methods, the resulting numerical solutions can exhibit unphysical values that can be another critical obstacle to accept such results by users. If a solution of the problem contains a steep gradient or even a shock or discontinuity, such unphysical oscillations in numerical solutions need not even diminish with a realistic refinement of discretization parameters as confirmed also in \cite{eimerImplicitFiniteVolume2022, barsukowImplicitActiveFlux2024}. Therefore,  to avoid such behavior, the higher order methods must be combined in general with lower order ones near critical positions in numerical solutions.

Such a combination of higher and lower order approximations is very well developed for explicit numerical methods involving many different numerical techniques. Among them, we mention the ones that we address later here, the Essentially Non-Oscillatory (ENO) approximations and their enhancement using WENO (Weighted ENO) approximations, see, e.g., \cite{shuEssentiallyNonoscillatoryWeighted1998}, and the high-resolution techniques involving different kind of limiters, see, e.g., \cite{levequeFiniteVolumeMethods2004}. 

An attractive approach to remove unphysical oscillations in numerical solutions is based on a posteriori techniques when, first, possibly oscillatory solutions with a high-order accurate method are computed and then, if undesirable features are observed in the numerical solution, lower-order methods are applied locally \cite{clainHighorderFiniteVolume2011,pimentel-garciaHighorderIncellDiscontinuous2024,maccaAlmostFailsafeAposteriori2024}. Unfortunately, in contrast to explicit methods, such an approach requires in the case of implicit methods some iterative procedure that is not easy to control to avoid all oscillations \cite{eimerImplicitFiniteVolume2022,michel-dansacTVDMOODSchemesBased2022,puppoQuinpiIntegratingStiff2024, maccaSemiimplicitTypeOrderAdaptiveCAT22025}.

When applying a priori techniques of variable accuracy for implicit discretization \cite{kuzminDesignGeneralpurposeFlux2006,duraisamyImplicitSchemeHyperbolic2007,kuzminLocallyBoundpreservingEnriched2020,arbogastThirdOrderImplicit2020,puppoQuinpiIntegratingConservation2022,quezadadelunaMaximumPrinciplePreserving2022}, they result in nonlinear algebraic systems that must be solved even for simple linear equations, e.g., the advection with constant speed. Moreover, to propose such approximation techniques that should keep the values of numerical solutions in physical bounds, it can not be done without involving downwind values, so the efficiency of fully upwinded numerical schemes for the usage of fast iterative solvers can be lost.

In \cite{frolkovicHighResolutionCompact2023}, a compact high-resolution implicit numerical scheme for one-dimen-sional hyperbolic systems is proposed that combines a second-order accurate WENO approximation in space with a first-order accurate limiter in time. The numerical solution is obtained by one forward and one backward substitution without iterations. 
Unfortunately, the time limiting towards the first accurate time discretization is strong for medium and large Courant numbers.
In this work, we relax this inefficiency by applying an inverse approach to the treatment of space and time variables. This technique in the form of an inverse Lax-Wendroff procedure is known and efficient when treating boundary conditions \cite{tanInverseLaxWendroffProcedure2010, shuInverseLaxWendroffBoundary2022}. In \cite{lackovaCompactSchemesAdvection2024}, it is used in combination with the direct scheme for one-dimensional advection. 

In the case of flow and transport problems on networks, when the flow enters in some input vertices and leaves it in several outflow vertices, such an inverse treatment of $x$ and $t$ can have several complementary advantages compared to standard direct treatment. 
At least for the constant speed advection equation, one can formally replace the roles of these two variables. Consequently, the treatment of initial conditions and the inflow (Dirichlet) boundary conditions can be replaced, and, to obtain the numerical solution, one can march in space from the inflow boundary toward the outflow. Using an analogous approach as in \cite{frolkovicHighResolutionCompact2023}, one can derive a compact ``inverse'' numerical scheme that can be applied to obtain the values of a numerical solution by a forward substitution.

For the inverse scheme, advanced approximation techniques with a larger stencil can be applied for time discretizations, and a simple two point stencil approximation can be used for the space approximation. To obtain higher accuracy with a smaller discretization stencil, we use uniform time steps, but space discretization can be variable. This significantly simplifies the treatment of boundary conditions when the larger stencils in space discretization can be difficult to apply for grid nodes near boundary with nonuniform grids. This is especially the case if a large gradient or even discontinuity of the solution enters the domain through boundary conditions \cite{shuInverseLaxWendroffBoundary2022}. In the case of advection on networks, this difficulty clearly occurs in the vertices that connect two edges \cite{eimerImplicitFiniteVolume2022,barsukowImplicitActiveFlux2024}.

The main advantages of this inverse approach that we exploit in this study is that the advection speed plays the role of ``slowness'', therefore large values of speed result in small values of slowness. Consequently, the additional ``time limiting'' known from direct schemes \cite{duraisamyConceptsApplicationTimeLimiters2003, duraisamyImplicitSchemeHyperbolic2007, arbogastThirdOrderImplicit2020,puppoQuinpiIntegratingConservation2022, frolkovicHighResolutionCompact2023,puppoQuinpiIntegratingStiff2024} is not necessary for large time steps, but only for ``too small''. 
We use this fact to propose an efficient space-time limiter for the inverse scheme that is based on the third-order accurate approximation in time and the second-order accurate approximation in space. 
This compact high resolution inverse scheme fulfills a priori and unconditionally the discrete minimum and maximum principle from which the stability follows automatically. 

Finally, to overcome the nonlinearity in algebraic equations introduced by the limiter in the case of a linear advection equation, we propose a simple predictor-corrector approach when both steps compute the value of numerical solution using an explicit formula. In theory, the bounds are preserved when the corrected value is identical to the predicted one, so in general more than one corrector step shall be applied; nevertheless, our limiter offers a priori a simple criterion when a correction preserves the bounds. In our experience, such an additional correction step for all chosen numerical examples had to be applied at most once and only occasionally, see related discussions in Section on numerical experiments.

We note that although linear advection-dominated problems have interesting applications for transport on networks \cite{borscheLocalTimeStepping2019, eimerImplicitFiniteVolume2022,sokacImpactSedimentLayer2021,veliskovaInverseTaskPollution2023}, our aim is to solve transport problems with nonlinear retardation coefficients \cite{kacurSolutionContaminantTransport2005,frolkovicSemianalyticalSolutionsContaminant2006,frolkovicNumericalSimulationContaminant2016,donatWENOSchemeCharacteristics2024}. Citing the work \cite{roeContributionsNumericalModelling1985} that ``almost any good numerical technique for the linear scalar wave equation could be converted, by a fairly simple mechanism, into an equally good numerical technique for a quasilinear system of conservation laws'', in this work we emphasize the linear case. However, we can confirm the quoted statement by presenting preliminary results for a nonlinear scalar example where the proposed numerical method exhibits similarly good behavior as in the linear case.

The paper is organized as follows. In Section \ref{sec-mathematical-model}, we introduce the basic notation of the network structure and the mathematical model. Section \ref{sec:numerical_schemes} is dedicated to the derivation of the proposed implicit numerical scheme, where we increase the accuracy of the scheme in sequential steps and discuss some important properties. Section \ref{sec:numerical_experiments} presents selected numerical experiments to illustrate the behavior and performance of the scheme. Finally, Section \ref{sec:conclusions} provides concluding remarks.

\section{Mathematical model}
\label{sec-mathematical-model}

We begin with the notation for networks that can represent in practice sewer systems \cite{markRiskAnalysesSewer1998,banik2014swmm5,veliskovaInverseTaskPollution2023} or district heating networks \cite{borscheLocalTimeStepping2019, mohring2021district,eimerImplicitFiniteVolume2022}. We suppose that a network is represented by a directed acyclic graph given by a set $\mathcal{V}$ of vertices and a set $\mathcal{E}$ of oriented edges with two vertices. The vertices $p_m \in \mathcal{V}$, $m=1,2,\ldots,M$ and the edges $l^e \in\mathcal{E}$, $e=1,2,\ldots,E$ are organized according to the following rules. 

Each edge $l^e$ of length $L^e$ is defined by a subset $\mathcal{V}^e \subset \mathcal{V}$ of two vertices, the left (input) vertex $p_{m_1} \in \mathcal{V}^e$ and the right (output) vertex $p_{m_2} \in \mathcal{V}^e$, such that $1 \le m_1 < m_2 \le M$. Consequently, there is no directed cycle in the graph. 

We distinguish the vertices $p_m$ in the network that are either only a left vertex or only a right vertex of some edge. 
Let $\mathcal{V}_{{in}}$ be a nonempty set containing the input vertices $p_m \in \mathcal{V}_{in}$ which are only the left points of some edges. Consequently, $p_1 \in \mathcal{V}_{{in}}$. Analogously, we denote the set $\mathcal{V}_{{out}}$ of vertices which are only the right points of some edges. Clearly, $p_M \in \mathcal{V}_{{out}}$. 

\begin{figure}[!ht]
    \centering
    \includegraphics[width=0.7\textwidth]{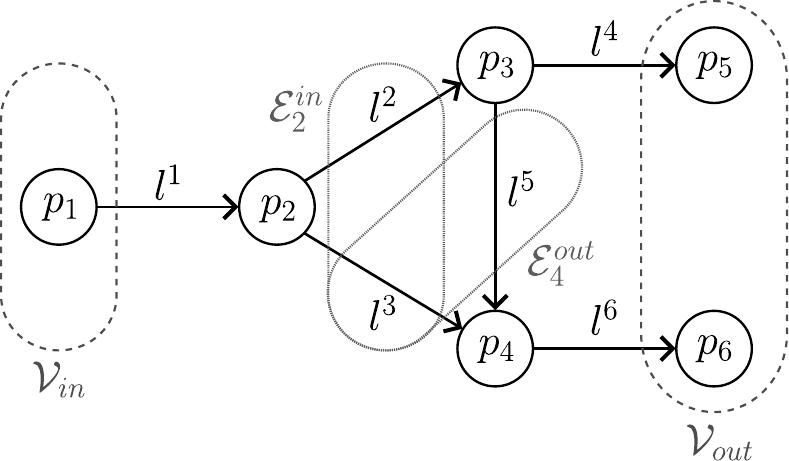}
    \caption{Schematic diagram of the network using the notation outlined in this section.}
    \label{fig:network_schema}
\end{figure}

Furthermore, by $\mathcal{E}_m^{in} \subset \mathcal{E}$, $p_m \in \mathcal{V}\setminus \mathcal{V}_{out}$, we denote nonempty subsets of edges $l^e \in \mathcal{E}_m^{in}$ for which $p_m$ is the left vertex, and $\mathcal{E}_m^{out} \subset \mathcal{E}$, $p_m \in \mathcal{V}\setminus \mathcal{V}^{in}$ are nonempty subsets of edges $l^e \in \mathcal{E}_m^{out}$ for which $p_m$ is the right vertex. There are no isolated vertices in the network. On the other hand, a vertex $p_m$ can be a left or right vertex of several edges.
Similarly to the numbering of vertices, we suppose that the edges are numbered such that if $l^{e_1} \in \mathcal{E}^{out}_m$ and $l^{e_2} \in \mathcal{E}^{in}_m$, then $e_1 < e_2$.
Consequently,  $l^1 \in \mathcal{E}_1^{in}$ and $l^E \in \mathcal{E}_M^{out}$. The structure of the network is illustrated in Figure \ref{fig:network_schema}.

Now we are ready to formulate the model equation to be solved on the network in the following form,
\begin{equation}
    \label{equation}
    \kappa^e \partial_t q^e + v^e(t) \partial_x q^e = 0 \,, \quad x \in (0,L^e) \,, \,\, t > 0 \,, \quad e=1,\ldots,E \,.
\end{equation}
Here, $q^e=q^e(x,t)$, $e=1,2,\ldots,E$ is the unknown function related to the edge $l^e$, where $t \in R$ is the time variable and $x \in [0,L^e]$ is the local space coordinate of the edge. The given flow rates $v^e>0$ can be variable with respect to time, and the capacity coefficient $\kappa^e$ \cite{levequeFiniteVolumeMethods2004}, for example, the cross-sectional area of the pipe, is a positive constant per each edge. In such a way, $q^e$ can represent the concentration of a transported quantity in the network, e.g. the mass per unit volume \cite{levequeFiniteVolumeMethods2004}.

For our purposes, we suppose that
\begin{equation}
    \label{init-condition}
    q^e(x,0)=0 \,, \quad x \in (0,L^e) \,, \quad e=1,2,\ldots, E \,,
\end{equation}
so initially there is no transported quantity in the network. Concerning coupling (or boundary) conditions, we proceed as follows.
For each inflow vertex $p_m \in \mathcal{V}_{in}$, we require a function $q_m(t)$ that defines the incoming concentration of the transported quantity into the network through the vertex $p_m$. For other vertices $p_m \in \mathcal{V}\setminus \mathcal{V}_{in}$, we define analogous functions $q_m(t)$ by summing the outcoming concentrations from all edges $l^e \in \mathcal{E}_m^{out}$, i.e.,
\begin{equation}
    \label{incomming-concentration}
    q_m(t) = \sum \limits_{l^e \in \mathcal{E}_m^{out}} q^e(L^e,t) \,.
\end{equation}
No boundary condition is required for $p_m \in \mathcal{V}_{out}$, but the function $q_m$ in \eqref{incomming-concentration} can be used to determine the outcoming concentration of the network.

For the vertices $p_m \in \mathcal{V}\setminus\mathcal{V}_{out}$, we assume the continuity of the concentration. 
In particular, we suppose that there are given non-negative coefficients $\alpha_m^{e}$, $l^e \in \mathcal{E}_m^{in}$ that sum to one with respect to $e$, and
\begin{equation}
    \label{coupling}
    q^{e}(0,t) = \alpha_m^{e} q_m(t) \,, l^{e} \in \mathcal{E}_m^{in} \,.
\end{equation}
If the vertex $p_m$ is the right vertex of a single edge $l^{e_1}$ and the left vertex of a single edge $l^{e_1}$ ($e_1<e_2$), then one has $\alpha_m^{e_1}=1$ and \eqref{incomming-concentration} and \eqref{coupling} turns to
\begin{eqnarray}
    \label{simple-coupling}
q^{e_2}(0,t)=q^{e_1}(L^{e_1},t) \,.
\end{eqnarray}

\begin{figure}[H]
  \centering
    \includegraphics[width=0.60\textwidth]{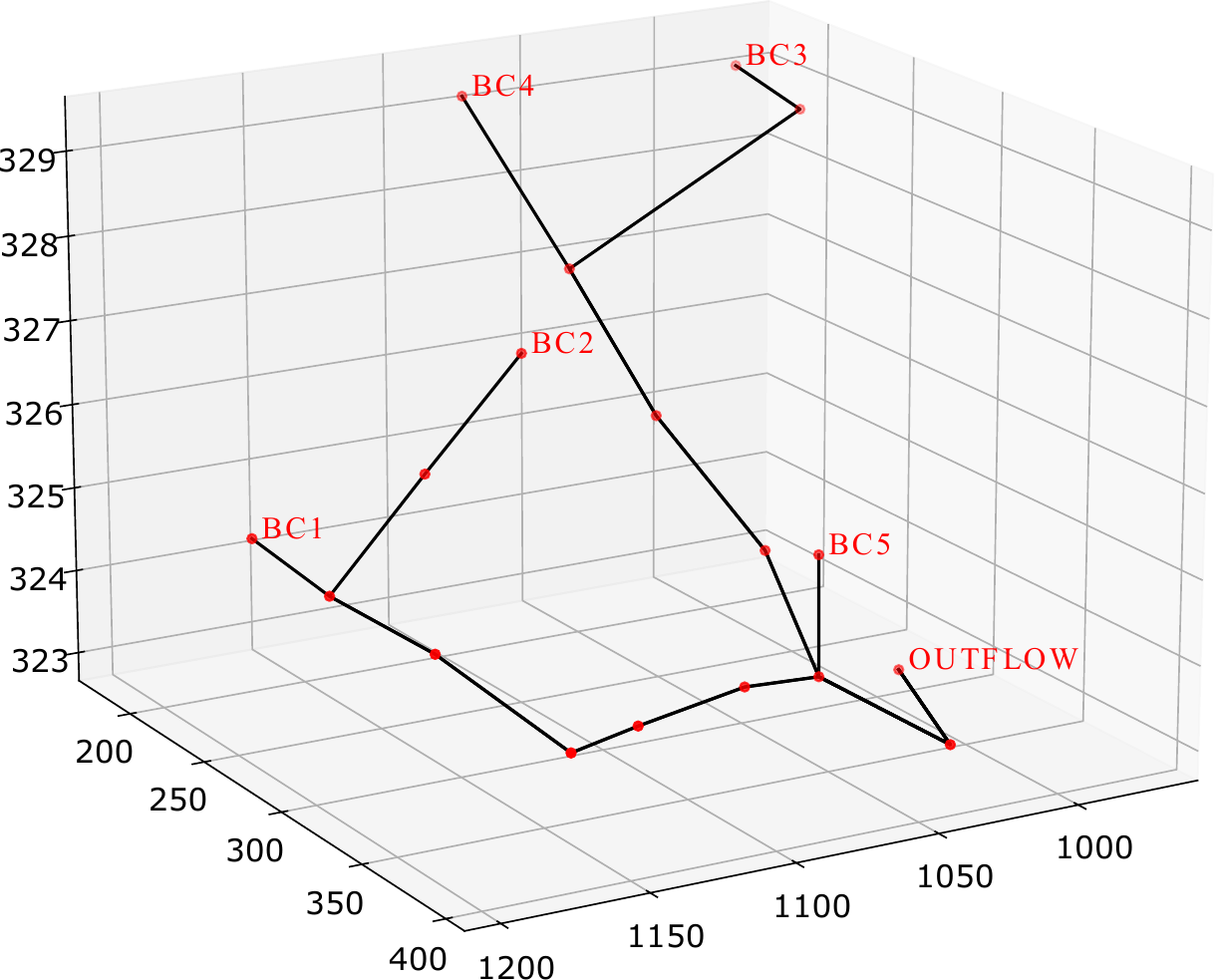}
  \caption{Three dimensional representation of a part of the sewer network in Revúca.}
  \label{fig:sewer3D}
\end{figure}

Our ultimate goal is to apply the proposed method to real-life examples of networks, such as the example of a sewer network in the Slovak city Revúca, see Figure \ref{fig:sewer3D}. Although it is described by a three-dimensional setup, one can approximate it well with a one-dimensional network in a plane where each edge represents a pipe with a different constant diameter. The flow regime can be estimated from given diameters and inclinations and then described by velocities assigned to pipes that can vary in time due to variable inflow into the network. Consequently, for such a flow regime one can model a pollutant transport driven by advection as described in this section. An example of related numerical simulations is given later in the section on numerical experiments.

\section{Numerical schemes}
\label{sec:numerical_schemes}

We begin with the description of our numerical scheme for the advection on a single edge $l^e$ of the network, therefore, we skip the index $e$ in the notation to simplify it,
\begin{equation}
    \label{advection-single-edge}
    \kappa \partial_t q + v(t) \partial_x q = 0 \,, \,\, x \in (0,L) \,, \,\, t > 0 \,. 
\end{equation}
The boundary and initial condition are given by
\begin{equation}
    \label{single-edge-conditions}
    q(0,t) = q_0(t) \,, \,\, t \ge 0 \,, \quad q(x,0) = 0 \,, \,\, x \in (0,L) \,.
\end{equation}

The equation \eqref{advection-single-edge} can be written in the form,
\begin{equation}
    \label{conservative-advection-single-edge}
    \partial_t \theta(q) + v(t) \partial_x q = 0
\end{equation}
with $\theta=\kappa q$ that we will follow in the notation for definitions of numerical schemes. As a byproduct, we can later consider nonlinear extensions of \eqref{conservative-advection-single-edge} where the coefficient $\theta$ can include some nonlinear dependence on the solution $q$. The purpose is to allow modeling of the influence of ``retardation'' of the transport as known, for example, for the transport with nonlinear sorption isotherms \cite{kacurSolutionContaminantTransport2005,frolkovicSemianalyticalSolutionsContaminant2006,frolkovicNumericalSimulationContaminant2016} or, eventually, due to immobile zones \cite{frolkovicNumericalSimulationContaminant2016,sokacImpactSedimentLayer2021}. As one can typically assume that $\theta'(q) \ge \kappa > 0$ for some given $\kappa$, formally, one can rewrite \eqref{conservative-advection-single-edge} into the form,
\begin{equation}
    \label{retardation-advection-single-edge}
    \partial_t q + \frac{v(t)}{\theta'(q)} \partial_x q = 0 \,,
\end{equation}
where the role of retardation is clearly demonstrated. Note that we do not consider such a nonlinear model here in a network, but only in the form \eqref{conservative-advection-single-edge}.

\subsection{Preliminaries}
\label{sec-preliminaries}

Firstly, we discretize the space and time variables. For $t \in [0,T]$ we use a uniform time step $\tau$ obtained for a chosen $N$ by $\tau=T/N$ and $t^n = n \tau$. For the space variable $x \in [0,L]$, we suppose a non-uniform grid given by $0=x_0 < x_1 < \ldots < x_I = L$ with $h_i = x_i - x_{i-1}$, $i=1,2,\ldots, I$. Our aim is to find approximations $Q_i^n \approx q(x_i,t^n)$ for $i=1,2,\ldots,I$ and $n=1,2,\ldots,N$ using the initial condition $Q_i^0=0$ and the boundary condition $Q_0^n=q_0(t^n)$.

To clearly introduce the method, we begin with the simplest implicit first-order accurate discretization of it. It can be derived, e.g., using the finite difference method,
\begin{equation}
    \label{first-order-scheme}
    \kappa \frac{Q_i^n - Q_i^{n-1}}{\tau} + v^n \frac{Q_i^n - Q_{i-1}^n}{h_i} = 0 \,,
\end{equation}
where $v^n:=v(t^n)>0$. It is useful to denote
the local Courant numbers,
\begin{equation}
    \label{courant-number}
    C_i^n = \frac{v^n \tau}{\kappa h_i} \,.
\end{equation}
when \eqref{first-order-scheme} can be written in the form
\begin{equation}
    \label{first-order-scheme-with-courant}
    Q_i^n - Q_i^{n-1} + C_i^n \left(Q_i^n - Q_{i-1}^n\right) = 0 \,.
\end{equation}
The standard application of the scheme \eqref{first-order-scheme-with-courant} is to march in time for $n=1,2,\ldots$ and using the given values of the initial condition $Q_i^0=0$ for $n=0$. 

We adopt a different approach. Namely, we initialize the values at the inflow boundary by
\begin{equation}
    \label{init-inflow-boundary}
    Q_0^n = q_0(t^n) \,, \,\, n=0,1,\ldots,N \,,
\end{equation}
and then apply the scheme \eqref{first-order-scheme} marching in space using $Q^0_i=0$ and the outer loop for the index $i$ with $i=1,2,\ldots,I$ to compute for each $i$
\begin{eqnarray}
    \label{marching-1st-order}
    Q_i^n = (1+C^n_i)^{-1} \left( Q_i^{n-1} + C^n_i Q_{i-1}^n \right) \,, \,\, n=1,2\ldots,N \,.
\end{eqnarray}

Next, we summarize the important properties of the first-order accurate method.
The numerical solution obtained with \eqref{marching-1st-order} fulfills the discrete minimum and maximum principle that we can write in the local form,
\begin{equation}
    \label{minmax-principle}
    \min \{ Q_i^{n-1},Q_{i-1}^n\} \le Q_i^n \le \max \{ Q_i^{n-1},Q_{i-1}^n\} \,.
\end{equation}
Consequently, the extremal values of the numerical solution can occur globally only at the initial values or the boundary values. 

The property \eqref{minmax-principle} is a straightforward consequence of \eqref{first-order-scheme-with-courant}, as it can be written in the ``positive coefficients'' form,
\begin{equation}
    \label{positive-coefficient-scheme}
    \alpha_i^n (Q_i^n - Q_i^{n-1}) + \beta_i^n (Q_i^n - Q_{i-1}^n) = 0 \,, \quad \alpha_i^n > 0 \,, \,\, \beta_i^n > 0 \,,
\end{equation}
with $\alpha_i^n=1$ and $\beta_i^n=C_i^n$ for \eqref{first-order-scheme-with-courant}.

We show later that high-resolution schemes can also be written in the form \eqref{positive-coefficient-scheme}, but with a nonlinear dependence of the coefficients $\alpha_i^n$ and $\beta_i^n$ on values $Q_i^n$ and $C_i^n$. Such schemes are called Local Extremum Diminishing (LED) in \cite{jamesonPositiveSchemesShock1995} as the maxima cannot increase and the minima cannot decrease in numerical solutions.

Having the notation and properties of the first-order accurate implicit finite difference scheme to solve \eqref{advection-single-edge} numerically, we can extend it to higher-order accuracy and high-resolution form. We do it in several steps.

\subsection{High order approximation in time}
\label{sec-high-order-time}

To easily introduce the extensions of the first order scheme \eqref{first-order-scheme}, we write the general scheme in the form where also the local conservation property is clearly given,
\begin{equation}
    \label{conservative-form}
    h_i \Theta_i^{n+1/2} + \tau v^n Q_i^n = h_i \Theta_i^{n-1/2} + \tau v^n Q_{i-1}^n \,.
\end{equation}
In the case of the first order accurate scheme \eqref{first-order-scheme}, one simply has
\begin{equation}
    \label{capacity-first-order}
    \Theta_i^{n+1/2} = \theta(Q_i^n) = \kappa Q_i^n  \,.
\end{equation}
For all numerical schemes presented here, the values $\Theta_i^{n+1/2}$ in \eqref{conservative-form} can be viewed as a representative concentration (e.g., averaged) for the time interval $(t^{n-1/2},t^{n+1/2})$ expressed in units of conserved quantity per unit length. Here, $t^{n+1/2} = t^n + \tau/2$. 

First, we aim to replace the first order approximation \eqref{capacity-first-order} by introducing a higher-order accurate approximation of $\partial_t q_i^{n}$ and keeping the first-order approximation in space. To obtain this accuracy, we extend the time discretization in \eqref{first-order-scheme} made of two values $Q_i^n$ and $Q_i^{n-1}$ with two additional unknowns, $Q_{i}^{n+1}$ and $Q_{i}^{n-2}$. 

In general, one has to avoid a fixed type of stencil for higher-order approximations that can lead to methods that produce numerical solutions with unphysical oscillations \cite{shuEssentiallyNonoscillatoryWeighted1998,levequeFiniteVolumeMethods2004}. Therefore, to approximate $\partial_t q_i^n$, we consider the following parametric approximations, 
\begin{equation}
    \label{parametric-time-approximation}
    \tau \partial_t q_i^n \approx Q_i^n - Q_{i}^{n-1} + \frac{1-w}{2} (Q_{i}^{n+1}- 2 Q_i^n + Q_{i}^{n-1}) + \frac{w}{2} (Q_i^n - 2 Q_{i}^{n-1} + Q_{i}^{n-2}) .
\end{equation}
This approximation is, in general, second-order accurate and reaches third-order accuracy if $w=1/3$ \cite{shuEssentiallyNonoscillatoryWeighted1998}. For two values, $w=0$ and $w=1$, it takes a reduced stencil, otherwise the four values, $Q_i^{n-2}$ up to $Q_i^{n+1}$, are used in the scheme. Note that later we consider a solution-dependent form of the parameters $w$ to obtain a high-resolution scheme. 

To write the scheme in the conservative form \eqref{conservative-form}, we replace \eqref{capacity-first-order} by
\begin{equation}
    \label{numerical-fluxes-2nd-order-in-time}
    \Theta_{i}^{n+1/2} = \kappa \left( Q_i^n + \frac{1-w}{2}(Q_{i}^{n+1} - Q_i^n) + \frac{w}{2} (Q_i^n - Q_{i}^{n-1})\right) \,.
\end{equation}
and in summary we obtain from \eqref{conservative-form}
\begin{eqnarray}
    \label{fully-implicit-inverse-scheme-fixed-w}
    Q_i^n + \frac{1-w}{2}(Q_{i}^{n+1} - Q_i^n) + \frac{w}{2} (Q_{i}^n - Q_{i}^{n-1}) + C_i^n Q_i^n = \\[1ex]
    \nonumber
    Q_i^{n-1} + \frac{1-w}{2}(Q_{i}^{n} - Q_i^{n-1}) + \frac{w}{2} (Q_{i}^{n-1} - Q_{i}^{n-2}) + C_i^n Q_{i-1}^n
    \,.
\end{eqnarray}

The system \eqref{fully-implicit-inverse-scheme-fixed-w} can be solved again by marching in space. In contrast to the first-order scheme, one has now to solve for each index $i$ in the outer loop a system of linear algebraic equations for unknowns $Q_i^1, Q_i^2, \ldots Q_i^N$ with a four-diagonal matrix if $w \in (0,1)$. For the first time step, when $n=1$, one also needs formally the values $Q_i^{-1}$ that can be obtained, for example, using the Lax-Wendroff procedure \cite{shuInverseLaxWendroffBoundary2022}. Note that in our case, as defined in \eqref{init-condition}, we simply assume that $q(x,t)=0$ for $t\le 0$.

Next, we aim to write the scheme in the form of \eqref{positive-coefficient-scheme}. To rewrite the scheme \eqref{fully-implicit-inverse-scheme-fixed-w} into that form, we use
\begin{equation}
    \label{numerical-flux-limiter-form-fully-implicit}
    \Theta_{i}^{n+1/2} = \kappa \left( Q_i^n + \frac{1}{2} \psi_i^{n} (Q_{i}^{n+1} - Q_i^n)\right) \,,
    \end{equation}
that is, if $Q_{i}^{n+1} \neq Q_i^n$, equivalent to \eqref{numerical-fluxes-2nd-order-in-time} by defining
\begin{equation}
    \label{limiter-to-parametric-for-fully-implicitm}
    \psi_i^{n} = 1 - w + w r_i^n \,, \quad r_i^n := \frac{Q_i^n-Q_{i}^{n-1}}{Q_{i}^{n+1}-Q_i^n} \,.
\end{equation}

Consequently, the scheme \eqref{fully-implicit-inverse-scheme-fixed-w} can be written in the form
\begin{equation}
    \label{scheme-first-order-time-second-space-coefficients}
    C_i^n \left(Q_i^n-Q_{i-1}^n\right) + \left( 1 + \frac{1}{2}\left(\frac{\psi_i^{n}}{r_i^n}-\psi_i^{n-1}\right)\right) (Q_i^n-Q_i^{n-1}) = 0 \,.
\end{equation}

It is now straightforward to formulate the conditions when the scheme \eqref{scheme-first-order-time-second-space-coefficients} is of the form \eqref{positive-coefficient-scheme} with nonnegative coefficients. Clearly, one has to require that 
\begin{equation}
    \label{first-conditon-coefficient-1t2x}
    -2 \le \frac{\psi_i^{n}}{r_i^n} - \psi_i^{n-1} \,.
\end{equation}
This condition is not satisfied for a fixed value of $w$ and for an arbitrary value of $r_i^n$. In fact, this can be seen as a source of violation of the discrete minimum and maximum principle in numerical solutions obtained with the scheme \eqref{fully-implicit-inverse-scheme-fixed-w}. At the same time, this is the motivation to introduce a solution-dependent form of the parameter $w=W(r_i^n)$ to follow WENO-type approximations \cite{shuEssentiallyNonoscillatoryWeighted1998} or to replace the values $\psi_i^{n}$ in \eqref{numerical-flux-limiter-form-fully-implicit} using the values of ``limiter function'' $\Psi(r_i^n)$ to follow the approach with limiters \cite{levequeFiniteVolumeMethods2004}.

We discuss advanced choices of the both approaches later for our compact inverse scheme, here we mention for \eqref{fully-implicit-inverse-scheme-fixed-w} the simplest choice of $w=W(r_i^n)$ that can be viewed as ENO (Essentially Non-Oscillatory) approximation,
\begin{equation}
    \label{ENO-1t2x}
    W(r_i^n) = \left \{\begin{array}{lr}
    0 \,, & |r_i^n| \ge 1 \\
    1 \,,  &  |r_i^n| \le 1
    \end{array}
    \right. 
    \quad \Rightarrow \quad
        \Psi(r_i^n) = \left \{\begin{array}{lr}
    1 \,, & |r_i^n| \ge 1 \\
    r_i^n \,,  &  |r_i^n| \le 1
    \end{array}
    \right. \,.
\end{equation}
One can verify by simple inspection that for \eqref{ENO-1t2x} the inequalities \eqref{first-conditon-coefficient-1t2x} are satisfied. Note that $\Psi(r_i^n) = 1 - W(r_i^n) +  W(r_i^n) r_i^n$.

The main difficulty of the scheme \eqref{scheme-first-order-time-second-space-coefficients} with \eqref{ENO-1t2x} is that it represents a nonlinear algebraic system of equations even in the case of the linear advection equation \eqref{advection-single-edge}. The system is fully coupled with the Jacobi matrix, which has a four-diagonal form. Of course, the scheme is still only first-order accurate in the space. In the next section, we improve not only the accuracy of the scheme, but also simplify its stencil in the implicit part.

\subsection{The inverse Lax-Wendroff procedure}
\label{sec-high-order-space}

To extend the scheme \eqref{scheme-first-order-time-second-space-coefficients} to higher-order accuracy in space, we apply the inverse Lax-Wendroff procedure in a partial form. In what follows, we use the abbreviation for the exact values: $q_i^n=q(x_i,t^n)$, $\partial_x q_i^n = \partial_x q(x_i,t^n)$, etc. The idea of the inverse Lax-Wendroff procedure is to replace the space derivatives in the Taylor series
\begin{equation}
\label{taylor-direct}
    v^n q_{i-1}^{n} = v^n q_i^n - h v^n \partial_x q_i^n + \frac{h^2}{2} v^n \partial_{xx} q_i^n + \mathcal{O}(h^3) \,, 
\end{equation}
by the derivatives of $q$ with respect to time. First, using the partial differential equation \eqref{conservative-advection-single-edge} itself,
\begin{equation}
\label{Lax-Wendroff-1}
  v^n \partial_x q_i^n = - \kappa \partial_t q_i^n
\end{equation}
and, second, by applying the space derivative to both sides,
\begin{equation}
\label{Lax-Wendroff-xx-to-tx}
  v^n \partial_{xx} q_i^n = - \kappa \partial_{tx} q_i^n \,.
\end{equation}

Before continuing with the derivation of the compact inverse scheme, we briefly discuss the standard inverse Lax-Wendroff procedure. The purpose is to show several advantages of our modification later. 

\begin{remark}
The fundamental concept of the inverse Lax-Wendroff procedure is to replace all space derivatives with time derivatives only, therefore one has to use the intermediate step by applying 
\begin{eqnarray}
\label{Lax-Wendroff-fully-implicit}
      - \kappa \partial_{tt} q_i^n = \partial_t \left(v^n \partial_{x} q_i^n\right) =& v^n \partial_{tx} q_i^n + \left(\frac{d}{dt} v^n\right)  \partial_{x} q_i^n  \\ \nonumber
      =& v^n \partial_{tx} q_i^n - \left(\frac{d}{dt} v^n\right) \frac{\kappa}{v^n} \partial_t q_i^n
\end{eqnarray}
and
\begin{eqnarray}
\label{Lax-Wendroff-fully-implicit-B}
  v^n \partial_{xx} q_i^n =  \frac{\kappa^2}{v^n} \partial_{tt} q_i^n - \left(\frac{\kappa}{v^n}\right)^2 \left(\frac{d}{dt} v^n\right) \partial_{t} q_i^n  \,.
\end{eqnarray}

Substituting all of the space derivatives with the time derivatives derived above, one obtains the following:
\begin{equation}
\label{taylor-direct-substitued}
    v^n q_{i-1}^{n} = v^n q_i^n + h \kappa \partial_t q_i^n - \frac{1}{2}\left(\frac{h \kappa}{v^n}\right)^2 \left(\frac{d}{dt} v^n\right) \partial_{t} q_i^n + \frac{(h \kappa)^2}{2 v^n} \partial_{tt} q_i^n  + \mathcal{O}(h^3) \,. 
\end{equation}

Clearly, if one considers the nonlinear equation \eqref{conservative-advection-single-edge}, the procedure becomes even more complicated. In any case, if the procedure is completed and previously introduced approximations are used, the stencil of the resulting numerical scheme is analogous to the one presented in Section \ref{sec-high-order-time} and one has still to solve a four-diagonal system of algebraic equations. Moreover, the scheme would involve derivatives of coefficients, so the conservative form \eqref{conservative-form} of the scheme can be lost.

\end{remark}

To obtain a compact inverse implicit numerical scheme, we keep the term containing the mixed derivatives $\partial_{tx} q_i^n$ in \eqref{taylor-direct-substitued} as used for the direct Lax-Wendroff procedure in several papers \cite{zorioApproximateLaxWendroff2017,baezaReprintApproximateTaylor2018, baezaApproximateImplicitTaylor2020,chouchoulisJacobianfreeImplicitMDRK2024,maccaSemiimplicitTypeOrderAdaptiveCAT22025}. Consequently, we obtain 
\begin{equation}
\label{taylor-compact-substituted}
    v^n q_{i-1}^{n} = v^n q_i^n + h \kappa \partial_t q_i^n - \frac{h^2}{2} \kappa \partial_{tx} q_i^n + \mathcal{O}(h^3) \,. 
\end{equation}

To approximate $\partial_t q_i^n$ in \eqref{taylor-compact-substituted}, we use the parametric discretization as given in \eqref{parametric-time-approximation}. To approximate $\partial_{tx} q_i^n$ in \eqref{taylor-compact-substituted}, it is enough to use a first-order accurate approximation that will not destroy the overall second-order accuracy of the derived scheme. Therefore, we choose
\begin{equation}
    \label{backward-space-finite-difference}
    h \partial_{tx} q_i^n \approx \partial_t q_i^n - \partial_t q_{i-1}^n
\end{equation}
and for $\partial_t q_i^n$ in \eqref{backward-space-finite-difference} we use
\begin{equation}
    \label{parametric-approximation-time-first-order}
    \tau \partial_t q_i^n \approx (1-w) (Q_{i}^{n+1} - Q_i^n) + w (Q_i^n - Q_{i}^{n-1}) \,,
\end{equation}
and analogously for $\partial_t q_{i-1}^n$. In summary, we obtain
 \begin{eqnarray}
    \label{summary}
\tau  \partial_t q_i^n - \tau \frac{h}{2} \partial_{tx} q_i^n \approx Q_i^n - Q_{i}^{n-1} + \\
\nonumber
 \frac{1-w}{2} \left(Q_{i-1}^{n+1} - Q_{i-1}^n -Q_i^n + Q_i^{n-1}\right) +  
 \frac{w}{2} \left(Q_{i-1}^{n} - Q_{i-1}^{n-1} - Q_i^{n-1} + Q_i^{n-2} \right) \,.
\end{eqnarray}

Consequently, one can express $\Theta_i^{n+1/2}$ in the conservative scheme \eqref{conservative-form} by
\begin{equation}
    \label{compact-fluxes-2nd-order-in-time}
    \Theta_{i}^{n+1/2} = \kappa \left( Q_i^n + \frac{1-w}{2}(Q_{i-1}^{n+1} - Q_i^n) + \frac{w}{2} (Q_{i-1}^n - Q_{i}^{n-1}) \right) \,.
\end{equation}
In summary, the scheme \eqref{conservative-form} with \eqref{compact-fluxes-2nd-order-in-time} can be written in the form
\begin{eqnarray}
    \label{compact-inverse-scheme-fixed-w}
    Q_i^n + \frac{1-w}{2}(Q_{i-1}^{n+1} - Q_i^n) + \frac{w}{2} (Q_{i-1}^n - Q_{i}^{n-1}) + C_i^n Q_i^n = \\[1ex]
    \nonumber
    Q_i^{n-1} + \frac{1-w}{2}(Q_{i-1}^{n} - Q_i^{n-1}) + \frac{w}{2} (Q_{i-1}^{n-1} - Q_{i}^{n-2}) + C_i^n Q_{i-1}^n
    \,.
\end{eqnarray}

The advantages of the scheme \eqref{compact-inverse-scheme-fixed-w} over the scheme \eqref{fully-implicit-inverse-scheme-fixed-w} are evident. First, the scheme \eqref{compact-inverse-scheme-fixed-w} gives the formal second-order accuracy not only in time but also in space. It preserves the parametric form of the approximation with reduced stencils for two choices of the parameter, $w=0$ and $w=1$ that we can use, see below, to define a WENO approximation. 
The most important difference is that the compact scheme \eqref{compact-inverse-scheme-fixed-w} does not contain the ``future'' value $Q_i^{n+1}$ for any value of $w$. Therefore, if the equations \eqref{compact-inverse-scheme-fixed-w} for any $i$ are solved sequentially for $n=1,2,\ldots,N$, they contain only a single unknown $Q_i^n$. Consequently, for the linear advection equation \eqref{advection-single-edge}, the value $Q_i^n$ can be expressed explicitly.

To propose a particular choice of $w \in [0,1]$, one can use $w=1/3$ which corresponds to the third order accurate approximation of the time derivative in \eqref{parametric-time-approximation}. In the case of constant coefficients $\theta =  \kappa q$ and $v^n \equiv v$ when also the Courant number is constant, i.e., $C_i^n \equiv C$, one can use $w = \bar w$ with
\begin{equation}
    \label{parameter-courant-number}
    \bar w = \frac{2 + \frac{1}{C}}{6} \,.
\end{equation}
The scheme \eqref{compact-inverse-scheme-fixed-w} with \eqref{parameter-courant-number} is of the third order accuracy for constant coefficients as proved for an analogous direct compact scheme in \cite{frolkovicSemiimplicitMethodsAdvection2022}. For too small Courant numbers $C< 1/4$, one can use $w=1$ instead of \eqref{parameter-courant-number}.

\begin{remark}
    \label{rem-weno}
Having the parametric form of the compact inverse implicit scheme \eqref{compact-inverse-scheme-fixed-w}, it is natural to ask if a straightforward application of a standard WENO approximation to define their solution dependent form can be used here. Namely, following \cite{shuEssentiallyNonoscillatoryWeighted1998}, one can consider the definitions
\begin{equation}
    \label{weno-parameters}
    w =  W(r) = \frac{\bar w}{\bar w + (1-\bar w) r} \,, \quad r = \frac{(\epsilon+(Q_{i-1}^{n+1}-Q_i^n)^2)^2}{(\epsilon+(Q_{i-1}^{n}-Q_i^{n-1})^2)^2} \,,
\end{equation}
where $\epsilon>0$ is a small parameter to be chosen. The behavior of $W(r)$ clearly follows the properties that $W(r)\rightarrow1$ for $r\rightarrow0$ and $W(r)\rightarrow0$ for $r\rightarrow\infty$ and $W(1)=\bar w$, where $r$ is analogous to a ``smoothness indicator'' \cite{shuEssentiallyNonoscillatoryWeighted1998} to choose appropriate weights between the extremal values $w=1$ and $w=0$ and the preferable choice $\bar w$. In \cite{shuEssentiallyNonoscillatoryWeighted1998,arbogastThirdOrderImplicit2020,puppoQuinpiIntegratingConservation2022,puppoQuinpiIntegratingStiff2024}, analogous indicators are based solely on a space reconstruction involving only several numerical values at time $t=t^n$, here the definition of $r$ in \eqref{weno-parameters} is rather heuristic.

We can report very good behavior of such WENO approximations to suppress small nonphysical negative values for the numerical solutions, but rather large clipping of local extremes. To illustrate this behavior, we use it and compare it later in a numerical example on a triangular network in Section \ref{sec:triangle}.
\end{remark}

As we mentioned in Remark \ref{rem-weno}, it is an open question how to propose the WENO approximation for \eqref{compact-inverse-scheme-fixed-w} with some appropriate smoothness indicators. In what follows, we prefer the variable accuracy approximation with A limiter function to derive a high-resolution form. Its derivation is well supported by some analysis, and the expected properties are confirmed later by chosen numerical experiments. 

\subsection{High-resolution compact inverse scheme}
\label{sec-high-resolution}

The core concept of the derivation of the high-resolution scheme is to rewrite \eqref{compact-inverse-scheme-fixed-w} in the form \eqref{positive-coefficient-scheme} with nonnegative coefficients. To obtain it, we proceed as follows.
Analogously to \eqref{numerical-flux-limiter-form-fully-implicit}, we rewrite \eqref{compact-fluxes-2nd-order-in-time} into the form
\begin{equation}
    \label{numerical-flux-limiter-form}
    \Theta_{i}^{n+1/2} = \kappa \left(Q_i^n + \frac{1}{2} \psi_i^{n} (Q_{i-1}^{n+1} - Q_i^n) \right) \,,
\end{equation}
that is, if $Q_{i-1}^{n+1} \neq Q_i^n$, equivalent to \eqref{compact-fluxes-2nd-order-in-time} by defining
\begin{equation}
    \label{limiter-to-parametric-form}
    \psi_i^{n} = 1 - w + w r_i^n \,, \quad r_i^n := \frac{Q_{i-1}^n-Q_{i}^{n-1}}{Q_{i-1}^{n+1}-Q_i^n} \,.
\end{equation}
Clearly, the indicator $r_i^n$ differs from the one defined in \eqref{limiter-to-parametric-for-fully-implicitm} for the fully implicit scheme \eqref{fully-implicit-inverse-scheme-fixed-w} and now depends only on the single unknown $Q_i^n$. 

Using \eqref{limiter-to-parametric-form}, the scheme \eqref{compact-inverse-scheme-fixed-w} can be written in the form
\begin{equation}
    \label{scheme-compact-with-psi}
    C_i^n \left(Q_i^n-Q_{i-1}^n\right) + Q_{i}^n-Q_i^{n-1} + \frac{1}{2}\left(\frac{\psi_i^{n}}{r_i^n}-\psi_i^{n-1}\right) (Q_{i-1}^n-Q_i^{n-1}) = 0 \,.
\end{equation}
The form \eqref{scheme-compact-with-psi} is nonlinear even for linear advection \eqref{advection-single-edge} due to the dependence of $r_i^n$ on $Q_i^n$. Note that the values $\psi_i^{n-1}$ are determined solely by the numerical values already computed when using the forward substitution.

Now, using 
$$
Q_{i-1}^n-Q_i^{n-1} = Q_i^n - Q_i^{n-1} - (Q_i^n - Q_{i-1}^n) \,,
$$
we obtain
\begin{eqnarray}
    \label{scheme-compact-with-coefficients}
    \left( C_i^n - \frac{1}{2}\left(\frac{\psi_i^{n}}{r_i^n}-\psi_i^{n-1}\right)\right) \left(Q_i^n-Q_{i-1}^n\right) + \\
    \nonumber\left( 1 + \frac{1}{2}\left(\frac{\psi_i^{n}}{r_i^n}-\psi_i^{n-1}\right)\right) (Q_{i}^n-Q_i^{n-1}) = 0 \,.
\end{eqnarray}

In the view of \eqref{scheme-compact-with-coefficients} it is straightforward to investigate the properties of the compact inverse implicit scheme \eqref{compact-inverse-scheme-fixed-w}. To obtain the nonnegative coefficients in \eqref{scheme-compact-with-coefficients}, one has to require
\begin{eqnarray}
    \label{inequalities}
    -2 \le \frac{\psi_i^{n}}{r_i^n}-\psi_i^{n-1} \le 2 C_i^n \,.
\end{eqnarray}
One can see that the left inequality in \eqref{inequalities} is analogous to the one in \eqref{first-conditon-coefficient-1t2x} for the fully implicit scheme \eqref{fully-implicit-inverse-scheme-fixed-w}. The bound in the right inequality in \eqref{inequalities} depends on the Courant numbers $C_i^n$, but in a favorable way, contrary to the bounds of the direct implicit schemes \cite{duraisamyImplicitSchemeHyperbolic2007,frolkovicHighResolutionCompact2023}, see Remark \ref{rem-weno}. Namely, if one can insure that $C_i^n \ge 1$, the inequalities are analogous to those used to construct many standard (Courant number independent) limiters in the literature \cite{levequeFiniteVolumeMethods2004,kemmComparativeStudyTVDlimiterswellknown2011,zhangReviewTVDSchemes2015}, 
\begin{eqnarray}
    \label{inequalities-reduced}
    -2 \le \frac{\psi_i^{n}}{r_i^n}-\psi_i^{n-1} \le 2 \,.
\end{eqnarray}
We later use this fact for the nonlinear extension of our scheme, because the Courant number in that case depends on the solution, so the limiters that are independent of Courant numbers are an attractive option. Nevertheless, for the linear case we will not use it to obtain the unconditionally local bounds preserving numerical scheme.

If $Q_i^n$ is the solution of the compact inverse scheme \eqref{compact-inverse-scheme-fixed-w} for some fixed value of $w$ and the inequalities \eqref{inequalities} are fulfilled, then the coefficients in \eqref{scheme-compact-with-coefficients} are nonnegative, and the discrete minimum and maximum principle \eqref{minmax-principle} is fulfilled. In fact, this is a simple indirect criterion (or a sufficient condition) for this property.
Of course, one cannot expect it in general, so for any fixed value of $w$ one can observe a violation of \eqref{minmax-principle} due to violation of \eqref{inequalities}.

To obtain a numerical scheme that preserves the bounds \eqref{minmax-principle} of numerical solutions, one can replace the values $\psi_i^{n}$ defined by \eqref{limiter-to-parametric-form} by the values of some limiter function ${\Psi}(r_i^n)$ which shall preserve the inequalities \eqref{inequalities}. 
The so called high-resolution form \cite{levequeFiniteVolumeMethods2004} of the compact inverse scheme \eqref{compact-inverse-scheme-fixed-w} is then obtained by using \eqref{conservative-form} with
\begin{equation}
    \label{theta-limiter-form}
    \Theta_{i}^{n+1/2} = \kappa \left( Q_i^n + \frac{{\Psi}(r_{i}^n)}{2}(Q_{i-1}^{n+1} - Q_i^n)\right) \,,
\end{equation}
or the form using a WENO type of approximation,
\begin{equation}
    \label{theta-weno-form}
    \Theta_{i}^{n+1/2} = \kappa \left( Q_i^n + \frac{1-W(r_i^n)}{2}(Q_{i-1}^{n+1} - Q_i^n) + \frac{W(r_i^n)}{2}(Q_{i-1}^{n} - Q_i^{n-1})\right) \,.
\end{equation}
These two definitions are equivalent if
\begin{equation}
    \label{psiin-from-W}
    \Psi(r_i^n) = 1 - W(r_i^n) + W(r_i^n) r_i^n \,,
\end{equation}
but they can give different results if some iterative (inexact) solver is used to solve the resulting nonlinear algebraic equation \eqref{conservative-form} with either \eqref{theta-limiter-form} or \eqref{theta-weno-form}.

Note that the standard ENO approximation \eqref{ENO-1t2x} of $\Psi$ or $W$ does not fulfill the Courant number dependent form \eqref{inequalities} for $C_i^n\ge 0$, but only the one \eqref{inequalities-reduced} independent of $C_i^n$. Inspired by \cite{aroraWellBehavedTVDLimiter1997,duraisamyImplicitSchemeHyperbolic2007,frolkovicHighResolutionCompact2023}, we propose a novel limiter function $\Psi=\Psi(r)$ for \eqref{compact-inverse-scheme-fixed-w} and the related WENO approximation $W=W(r)$.

To derive it, we require the condition that is very often used to construct limiter functions,
\begin{equation}
    \label{1st-inequality}
    0 \le \Psi(r) \le 2 \,, \quad r \in R \,.
\end{equation}
To construct the limiter function for an available value of $\Psi(r_i^{n-1})$ to fulfill \eqref{inequalities}, we require
\begin{eqnarray}
    \label{2nd-inequality}
    0 \le \frac{\Psi(r)}{r} \le 2 C_i^n + \Psi(r_i^{n-1})\,.
\end{eqnarray}
These two requirements, \eqref{1st-inequality} and \eqref{2nd-inequality}, together with the preferable form \eqref{limiter-to-parametric-form} for some chosen linear weight $\bar w$, motivates the following definition,
\begin{equation}
    \label{our-limiter}
    \Psi(r) = \max\{0, \min\{(2 C_i^n + \Psi(r_i^{n-1}) r, 1-\bar w + \bar w r, 2 \}\} \,.
\end{equation}
For accuracy reasons \cite{levequeFiniteVolumeMethods2004,kemmComparativeStudyTVDlimiterswellknown2011,zhangReviewTVDSchemes2015}, one expects that if $r=1$ then $\Psi(1)=1$. If $\Psi(r_i^{n-1}) = 0$ and $C_i^n < 1/2$, the right inequality in \eqref{inequalities} requires that $\Psi(r)<1$, so this accuracy requirement need not to be met and the formal second order accuracy could be locally decreased towards the first order accurate scheme when $\Psi(r_i^n)=0$. Clearly, this can happen only for too small Courant numbers.

The final scheme takes then the form,
\begin{eqnarray}
    \label{compact-scheme-with-limiters}
    Q_i^n + \frac{1}{2}\Psi(r_i^n)(Q_{i-1}^{n+1}-Q_i^n) +C_i^n Q_i^n = \\ \nonumber
    Q_i^{n-1} + \frac{1}{2}\Psi(r_i^{n-1})(Q_{i-1}^{n}-Q_i^{n-1}) +C_i^n Q_{i-1}^n \,.
\end{eqnarray}
The scheme \eqref{compact-scheme-with-limiters} can be obtained using \eqref{psiin-from-W} by the definition of the solution dependent form of the parameter $w$,
\begin{equation}
    \label{psi-from-weno}
    W(r) = \left\{
\begin{array}{lr}
\frac{1}{r-1}     & r \ge \frac{1+\bar w}{\bar w}  \\[1ex]
 \bar w    &  r \ge \frac{1-\bar w}{2 C_i^n + \Psi(r_i^{n-1}) - \bar w} \\[1ex]
 \frac{1-(2 C_i^n + \Psi(r_i^{n-1})) r}{1 - r} &  r \ge 0 \\[1ex]
 \frac{1}{1-r} & r \le 0 
\end{array}
    \right.   \,.
\end{equation}

To solve the nonlinear algebraic equation \eqref{compact-scheme-with-limiters}, we propose the following predictor-corrector approach. Firstly, we compute a predicted value $\tilde Q_i^n$ by using the preferred value of $w=\bar w$, namely by solving the linear alebraic equation for the single unknown $Q_i^n$, 
\begin{eqnarray}
    \label{predictor}
    Q_i^n + \frac{1-\bar w}{2}(Q_{i-1}^{n+1}-Q_i^n) + \frac{\bar w}{2}(Q_{i-1}^{n}-Q_i^{n-1}) + C_i^n Q_i^n = \\ \nonumber
    Q_i^{n-1} + \frac{1}{2} \Psi(\tilde r_i^{n-1})(Q_{i-1}^{n}-Q_i^{n-1}) + C_i^n Q_{i-1}^n \,.
\end{eqnarray}
Note that to ensure the conservative form of our scheme, we have taken into account that the previous values of $\Theta_i^{n-1/2}$ were computed using the value $\Psi(\tilde r_i^{n-1})$, see below. It is important that $\Psi(\tilde r_i^{n-1})$ fulfills the inequalities \eqref{1st-inequality}.

Having $\tilde Q_i^n$ from \eqref{predictor}, we compute the value $\tilde r_i^n$ and $\Psi(\tilde r_i^n)$ by
\begin{eqnarray}
    \label{our-limiter-predicted}
    \tilde r_i^n = \frac{Q_{i-1}^{n}-Q_i^{n-1}}{Q_{i-1}^{n+1}-\tilde Q_i^{n}} \,, \\ \nonumber
    \Psi(\tilde r_i^{n}) = \max\{0,\min\{(2C_i^n+\Psi(\tilde r_i^{n-1})) \tilde r_i^n, 1-\bar w+\bar w\tilde r_i^n, 2\}\} \,.
\end{eqnarray}
Note that $0\le \Psi(\tilde r_i^n) \le 2$. Clearly, if one gets from \eqref{our-limiter-predicted} that $\Psi(\tilde r_i^n) = 1-\bar w+\bar w\tilde r_i^n$, the predicted value can be accepted because the correction below will return the same value. If not, the corrected value is computed by \eqref{conservative-form} using \eqref{theta-limiter-form}, it means from the linear algebraic equation for the single unknown $Q_i^n$,
\begin{eqnarray}
    \label{corrector}
    Q_i^n + \frac{1}{2}\Psi(\tilde r_i^n)(Q_{i-1}^{n+1}-Q_i^n) + C_i^n Q_i^n = \\ \nonumber
    Q_i^{n-1} + \frac{1}{2}\Psi(\tilde r_i^{n-1})(Q_{i-1}^{n}-Q_i^{n-1}) + C_i^n Q_{i-1}^n \,.
\end{eqnarray}

As \eqref{corrector} is only an approximation of the exact nonlinear form \eqref{compact-scheme-with-limiters}, in theory another correction might be necessary in general. It can be fairly easily decided by recomputing the value of $r_i^n$ using the corrected value $Q_i^n$ and checking the inequalities \eqref{2nd-inequality}. If they are not fulfilled, another corrector step can be computed. For the chosen examples, as we discuss in the section on numerical experiments, at most a single repetition of the corrector step is occasionally necessary.

\begin{remark}
    \label{rem-comparison}
To clearly show later in numerical experiments the advantages of the compact high-resolution inverse scheme \eqref{corrector} over the analogous direct scheme \cite{frolkovicHighResolutionCompact2023}, we present it here for completeness for the case of a constant Courant number $C\equiv C_i^n$. The predictor is computed by
\begin{eqnarray}
    \label{predictor-direct}
    Q_i^n + C \left(Q_i^n + \frac{1-\bar w}{2}(Q_{i+1}^{n-1}-Q_i^{n}) + \frac{\bar w}{2}(Q_{i}^{n-1}-Q_{i-1}^{n})\right) = \\ \nonumber
    Q_i^{n-1} + C \left(Q_{i-1}^n + \frac{\Phi(\tilde s_{i-1}^n)}{2}(Q_{i}^{n-1}-Q_{i-1}^{n}) \right) \,.
\end{eqnarray}
Afterwards, we compute 
\begin{equation}
    \label{our-limiter-direct}
    \tilde s_i^n = \frac{Q_{i}^{n-1}-Q_{i-1}^{n}}{Q_{i+1}^{n-1}-\tilde Q_i^{n}} \,, \quad 
    \Phi(\tilde s_i^{n}) = \max\{0,\min\{(\frac{2}{C}+\Phi(\tilde s_{i-1}^{n})) \tilde s_i^n, 1-\bar w+\bar w\tilde s_i^n, 2\}\}
\end{equation}
Finally, the corrected value is computed by
\begin{eqnarray}
    \label{corrector-direct}
    Q_i^n + C \left( Q_i^n + \frac{\Phi(\tilde s_i^n)}{2}(Q_{i+1}^{n-1}-Q_i^n) \right) = 
    \\ \nonumber Q_i^{n-1} + C \left(Q_{i-1}^n + \frac{\Phi(\tilde s_{i-1}^n)}{2}(Q_{i}^{n-1}-Q_{i-1}^{n}) \right) \,.
\end{eqnarray}
Notice that the Courant number $C$ in the definition \eqref{our-limiter-direct} of $\Phi$ is in a denominator, which is unfavorable situation for the limiter function $\Psi$ in the case of large Courant numbers. For $C\ge 4$, the preferable choice for the parameter $\bar w$ is $1$.
\end{remark}

\subsection{The nonlinear case}
\label{sec-nonlinear}

Following the notation of the advection equation given in \eqref{conservative-advection-single-edge}, it is rather straightforward to extend the compact inverse scheme for the case of $\theta=\theta(q)$ including a nonlinear dependence. For that purpose, we denote $\Theta_i^n := \theta(Q_i^n)$ and one has to express $\Theta_i^{n+1/2}$ in the conservative scheme \eqref{conservative-form} by
\begin{equation}
    \label{compact-fluxes-2nd-order-in-time-nonlinear}
    \Theta_{i}^{n+1/2} =  \Theta_i^n + \frac{1-w}{2}(\Theta_{i-1}^{n+1} - \Theta_i^n) + \frac{w}{2} (\Theta_{i-1}^n - \Theta_{i}^{n-1}) \,.
\end{equation}
Consequently, the scheme \eqref{conservative-form} with \eqref{compact-fluxes-2nd-order-in-time} can be written in the form
\begin{eqnarray}
    \label{compact-inverse-scheme-fixed-w-nonlinear}
    \Theta_i^n + \frac{1-w}{2}(\Theta_{i-1}^{n+1} - \Theta_i^n) + \frac{w}{2} (\Theta_{i-1}^n - \Theta_{i}^{n-1}) + v^n Q_i^n = \\[1ex]
    \nonumber
    \Theta_i^{n-1} + \frac{1-w}{2}(\Theta_{i-1}^{n} - \Theta_i^{n-1}) + \frac{w}{2} (\Theta_{i-1}^{n-1} - \Theta_{i}^{n-2}) + v^n Q_{i-1}^n
    \,.
\end{eqnarray}

Analogously to \eqref{numerical-flux-limiter-form}, we rewrite \eqref{compact-fluxes-2nd-order-in-time-nonlinear} into the form
\begin{equation}
    \label{numerical-flux-limiter-form-nonlinear}
    \Theta_{i}^{n+1/2} = \Theta_i^n + \frac{1}{2} \psi_i^{n} (\Theta_{i-1}^{n+1} - \Theta_i^n) \,,
\end{equation}
that is, if $Q_{i-1}^{n+1} \neq Q_i^n$, equivalent to \eqref{compact-fluxes-2nd-order-in-time-nonlinear} by defining
\begin{equation}
    \label{limiter-to-parametric-form-nonlinear}
    \psi_i^{n} = 1 - w + w r_i^n \,, \quad r_i^n := \frac{\Theta_{i-1}^n-\Theta_{i}^{n-1}}{\Theta_{i-1}^{n+1}-\Theta_i^n} \,.
\end{equation}

Using \eqref{limiter-to-parametric-form-nonlinear}, the scheme \eqref{compact-inverse-scheme-fixed-w-nonlinear} can be written in the form
\begin{equation}
    \label{scheme-compact-with-psi-nonlinear}
    \frac{\tau v^n}{h_i} \left(Q_i^n-Q_{i-1}^n\right) + \Theta_{i}^n-\Theta_i^{n-1} + \frac{1}{2}\left(\frac{\psi_i^{n}}{r_i^n}-\psi_i^{n-1}\right) (\Theta_{i-1}^n-\Theta_i^{n-1}) = 0 \,.
\end{equation}

Using, if $Q_i^n \neq Q_{i-1}^n$,
$$
Q_i^n - Q_{i-1}^n = \frac{Q_i^n - Q_{i-1}^n}{\Theta_i^n - \Theta_{i-1}^n}(\Theta_i^n - \Theta_{i-1}^n) \,,
$$
and
$$
\Theta_{i-1}^n-\Theta_i^{n-1} = \Theta_i^n - \Theta_i^{n-1} - \left(\Theta_i^n - \Theta_{i-1}^n \right) \,,
$$
we obtain
\begin{eqnarray}
    \label{scheme-compact-with-coefficients-nonlinear}
    \left( C_i^n - \frac{1}{2}\left(\frac{\psi_i^{n}}{r_i^n}-\psi_i^{n-1}\right)\right) \left(\Theta_i^n-\Theta_{i-1}^n\right) + \\
     \nonumber\left( 1 + \frac{1}{2}\left(\frac{\psi_i^{n}}{r_i^n}-\psi_i^{n-1}\right)\right) (\Theta_{i}^n-\Theta_i^{n-1}) = 0 \,,
\end{eqnarray}
where we have to extend the definition of Courant number to the nonlinear form,
\begin{equation}
    \label{Courant-nonlinear}
    C_i^n := \left \{
    \begin{array}{lr}
       \frac{Q_i^n - Q_{i-1}^n}{\Theta_i^n - \Theta_{i-1}^n} \frac{\tau v^n}{h_i}  &  Q_i^n \neq Q_{i-1}^n\\[1ex]
      \frac{\tau v^n}{\theta'(Q_i^n) h_i} & Q_i^n = Q_{i-1}^n 
    \end{array}
    \right.
     \,,
\end{equation}
which is compatible with \eqref{courant-number} for the linear case $\theta(q)=\kappa q$. Clearly, we have obtained the scheme of the form \eqref{positive-coefficient-scheme} with the coefficients that can be made nonnegative if the limiter from the previous section is used. Note that to make the limiter less computationally expensive and more robust, one can replace the Courant numbers $C_i^n$ in \eqref{scheme-compact-with-coefficients-nonlinear} with a minimal one $C^{\min}$ if such an estimate is available.

\section{Numerical experiments}
\label{sec:numerical_experiments}

The following section is dedicated to numerical experiments that aim to demonstrate different properties of numerical methods described in this paper. 

We start by evaluating the accuracy of the numerical schemes using a simple example with constant velocity and smooth solution in Section \ref{sec-eoc}. In Section \ref{sec-nonsmooth}, we test a well-known benchmark to assess the oscillation resilience of the methods not only for linear advection, but also for the nonlinear transport equation. Section \ref{sec:triangle} follows with an accuracy test using a nontrivial example taken from \cite{barsukowImplicitActiveFlux2024} computed on a triangular network, and finally, we apply the methods inspired by a real-world case involving a sewer network in the city of Revúca in Section \ref{sec-revuca}.

\subsection{Convergence study for a smooth solution}
\label{sec-eoc}

We begin with the simplest example of advection \eqref{advection-single-edge} having constant unit velocity with the exact solution $q(x,t)=q^0(x-t)$, where $q^0$ is some chosen function. 
Firstly, 
\begin{equation}
\label{init-exp}
    q^0(x)=\exp(-40(x+1)^2) \,, \,\, x \in R \,,
\end{equation}
and we solve the equation for $x \in (0,2)$ and $t \in [0,2]$. The boundary condition at $x=0$ is prescribed by the exact solution. In Table \ref{tab-01-exp-eoc} we present the results for the scheme \eqref{compact-inverse-scheme-fixed-w} with two different fixed values of the parameter $w$, $\bar w=1/3$ (the second order scheme) and $\bar w$ in \eqref{parameter-courant-number} (the third order scheme). Furthermore, the results of the high-resolution scheme \eqref{compact-scheme-with-limiters} with one corrector \eqref{corrector} are given. We present results for two Courant numbers to illustrate the stability and the EOC (Experimental Order of Convergence) for this example with the erro computed by
\begin{equation}
    \label{error}
    E = E_I^N = \sum \limits_{i,n} |Q_i^n - q(x_i,t^n)| \,.
\end{equation} 

\begin{figure}[!ht]
  \centering 
    \subfloat[$I=64, C=1$]{\includegraphics[width=0.49\textwidth]{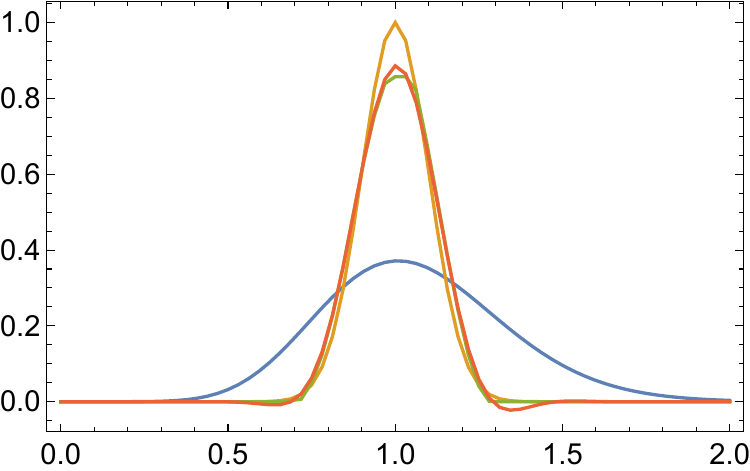}\label{fig:1}}
    \subfloat[$I=256, C=8$]{\includegraphics[width=0.49\textwidth]{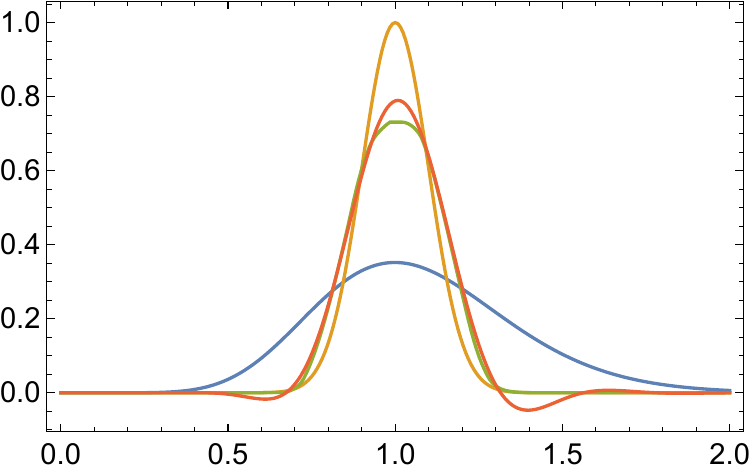}\label{fig:2}}
  \caption{Comparison of the example in Section \ref{sec-eoc}: the four curves correspond to the (interpolated) exact solution (orange), the 3rd order method (red), the high-resolution method (green), and the first-order accurate method (blue). On the left, the solution is obtained on a coarse grid with a small Courant number, while on the right, it is computed on a fine mesh with a large Courant number.}
  \label{fig:smooth}
\end{figure}

\begin{figure}[!ht]
  \centering 
    \subfloat[Unlimited numerical solutions]{\includegraphics[width=0.49\textwidth]{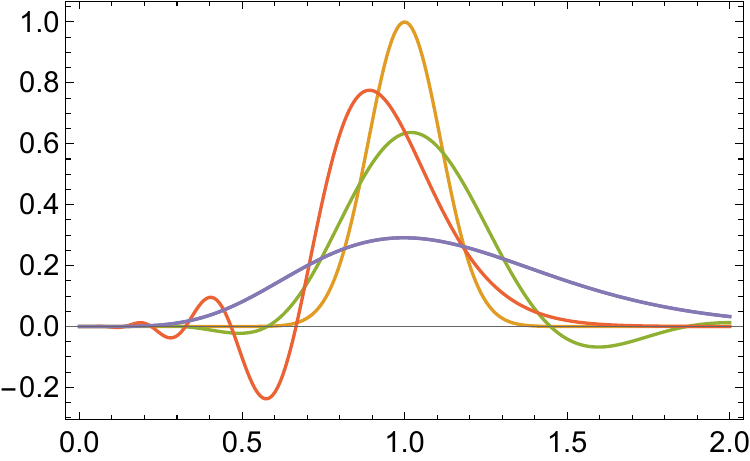}\label{fig:3}}
    \subfloat[Limited numerical solutions]{\includegraphics[width=0.49\textwidth]{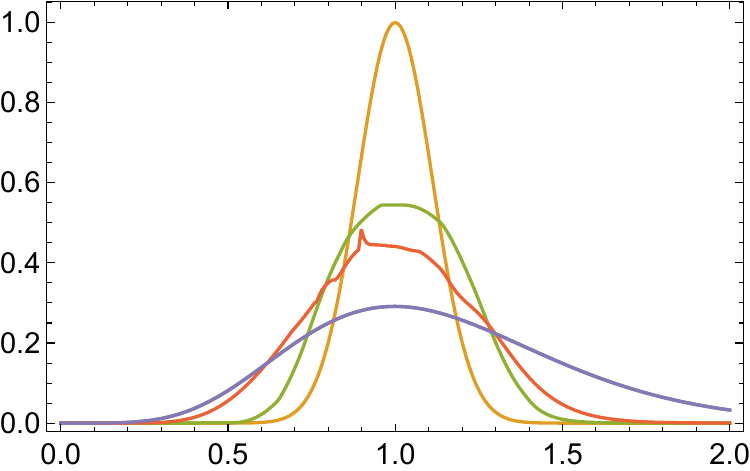}\label{fig:4}}
  \caption{Comparisons for the example in Section \ref{sec-eoc}: the curves correspond to the exact solution (orange), the direct method in Remark \ref{rem-comparison} (red), the inverse method (green), and the first-order method (blue), computed for $I=256$ and the Courant number of 16. The left plot shows the solution with a fixed parameter (unlimited), while the right one the high-resolution (limited) form.}
  \label{fig:smooth2}
\end{figure}

A visual comparison of the numerical solutions given by various methods is shown in Figures \ref{fig:smooth} and \ref{fig:smooth2}. 
One can see that the expected EOCs are obtained for each scheme. For illustration, the results of the first order scheme \eqref{first-order-scheme} are also given, which are significantly less precise as can be seen in Table \ref{tab-01-exp-eoc} and in Figures \ref{fig:smooth} and \ref{fig:smooth2}. 
Although the third-order method provides the best accuracy, the high-resolution method is more effective in suppressing negative oscillations, as clearly documented in Table \ref{tab-02-exp-min}.

\begin{table}[!ht]
\small
    \begin{center}
    \begin{tabular}{c c | c c | c c | c c | c c }
    \hline 
    \multicolumn{2}{c}{} & \multicolumn{8}{c}{$C = 2$}  \\
    \hline
    \multicolumn{2}{c}{} & \multicolumn{2}{c}{1st order} & \multicolumn{2}{c}{2nd order} & \multicolumn{2}{c}{3rd order} & \multicolumn{2}{c}{HR}  \\
    \hline
    $I $ & $N$  & E & EOC  & E & EOC  & E  & EOC & E  & EOC   \\ 
    \hline
    256  & 128  & 86.8 & - & 2.879  & - & 1.491  & -   & 1.573 & -   \\    
    512  & 256  & 51.3 & 0.76 & 0.639  & 2.17 & 0.191  & 2.96 &  0.198 & 2.99    \\
    1024  & 512  & 28.6  & 0.85 & 0.153  & 2.06 & 0.024 & 2.99   & 0.031 & 2.70   \\
    2048  & 1024  & 15.2  & 0.91 & 0.038  & 2.02 & 0.003  & 3.00  & 0.006 & 2.41  \\
    \hline 
    \hline 
    \multicolumn{2}{c}{} & \multicolumn{8}{c}{$C = 8$}  \\
    \hline
    \multicolumn{2}{c}{} & \multicolumn{2}{c}{1st order} & \multicolumn{2}{c}{2nd order} & \multicolumn{2}{c}{3rd order}  & \multicolumn{2}{c}{HR} \\
    \hline
    $I $ & $N$  & E & EOC  & E & EOC  & E  & EOC  & E  & EOC \\ 
    \hline
    256  & 32  & 176. & - & 41.86  & - & 40.23  & -  & 38.07 & -  \\    
    512  & 64  & 115. & 0.61 & 7.418  & 2.50 & 7.061  & 2.51  &  8.258  & 2.21  \\
    1024  & 128  & 70.8  & 0.71 & 1.089  & 2.77 & 0.963 & 2.87  & 1.105 & 2.90 \\
    2048  & 256  & 40.6  & 0.80 & 0.171  & 2.67 & 0.122  & 2.98  & 0.152 & 2.86     \\
    \hline 
    \end{tabular}
    \end{center}
    \caption{Comparison of the numerical errors $E \cdot 10^{-3}$ and the EOCs for the first, second, third order accurate schemes and the high resolution (HR) scheme for the constant speed and the initial function \eqref{init-exp} using different Courant numbers.}
    \label{tab-01-exp-eoc}
\end{table}

\begin{table}[!ht]
    \begin{center}
    \begin{tabular}{c c | c c | c c | c c}
    \hline 
    \multicolumn{2}{c}{} & \multicolumn{6}{c}{$C = 2$}  \\
    \hline
    \multicolumn{2}{c}{}  & \multicolumn{2}{c}{2nd order} & \multicolumn{2}{c}{3rd order} & \multicolumn{2}{c}{HR}  \\
    \hline
    $I $ & $N$  & $\min$ & $\max$  & $\min$ & $\max$  & $\min$  & $\max$  \\ 
    \hline
    256  & 128  & -8.1$\cdot10^{-5}$ & 0.9907 & -6.8$\cdot10^{-6}$  & 0.9896 & 7$\cdot10^{-17}$  & 0.979     \\    
    512  & 256  & -5.2$\cdot10^{-9}$ & 0.9988 & -2.6$\cdot10^{-8}$  & 0.9986 & 7$\cdot10^{-18}$ & 0.9951     \\
    1024  & 512  & -1.7$\cdot10^{-7}$  & 0.9999 & -4.6$\cdot10^{-8}$  & 0.9998 & 5$\cdot10^{-18}$ & 0.9985   \\
    2048  & 1024  & -5.1$\cdot10^{-8}$  & 1.0000 & -5.9$\cdot10^{-8}$  & 1.0000 & 5$\cdot10^{-18}$ & 0.9995    \\
    \hline 
    \hline 
    \multicolumn{2}{c}{} & \multicolumn{6}{c}{$C = 8$}  \\
    \hline
    \multicolumn{2}{c}{} & \multicolumn{2}{c}{2nd order}  & \multicolumn{2}{c}{3rd order} & \multicolumn{2}{c}{HR}  \\
    \hline
    $I $ & $N$  & $\min$ & $\max$  & $\min$ & $\max$  & $\min$  & $\max$  \\ 
    \hline
    256  & 32  & -4.7$\cdot 10^{-2}$ & 0.836 & -3.9$\cdot 10^{-2}$  & 0.832 & 3$\cdot 10^{-14}$  & 0.775   \\    
    512  & 64  & -4.7$\cdot 10^{-3}$ & 0.959 & -3.4$\cdot 10^{-3}$  & 0.957 & 5$\cdot 10^{-17}$  & 0.927   \\
    1024  & 128  & -3.7$\cdot 10^{-7}$  & 0.994 & -1.5$\cdot 10^{-7}$  & 0.993 & 6$\cdot 10^{-18}$ & 0.981   \\
    2048  & 256  & -1.8$\cdot 10^{-7}$  & 0.999 & -4.2$\cdot 10^{-7}$ & 0.999 & 5$\cdot 10^{-18}$ & 0.995      \\
    \hline 
    \end{tabular}
    \end{center}
    \caption{Comparison of the minimal and maximal values of numerical solutions at $t=2$ for the second and third order accurate schemes and the high resolution (HR) scheme for the constant velocity and the initial function \eqref{init-exp} using different Courant numbers.}
    \label{tab-02-exp-min}
\end{table}

Finally, we illustrate with this example a better behavior of the compact inverse schemes with respect to the compact direct scheme \cite{frolkovicHighResolutionCompact2023} as presented in Remark \ref{rem-comparison}. For that purpose, we compute the example with large Courant number $C_i^n \equiv C =16$. As can be seen in Figure \ref{fig:smooth2}, for the fixed parameter case, the nonphysical oscillations are much larger for the direct scheme than for the inverse scheme. The high-resolution (HR) version of the direct scheme has difficulties to capture well the large gradient of the numerical solution that enters the computational interval through the inflow boundary for large Courant numbers. Moreover, the numerical solution obtained with the HR direct scheme is more limited towards the first order scheme than the HR inverse scheme.

\subsection{Benchmark with a non-smooth solution}
\label{sec-nonsmooth}
\begin{figure}[!ht]
  \centering
   {\includegraphics[width=0.49\textwidth]{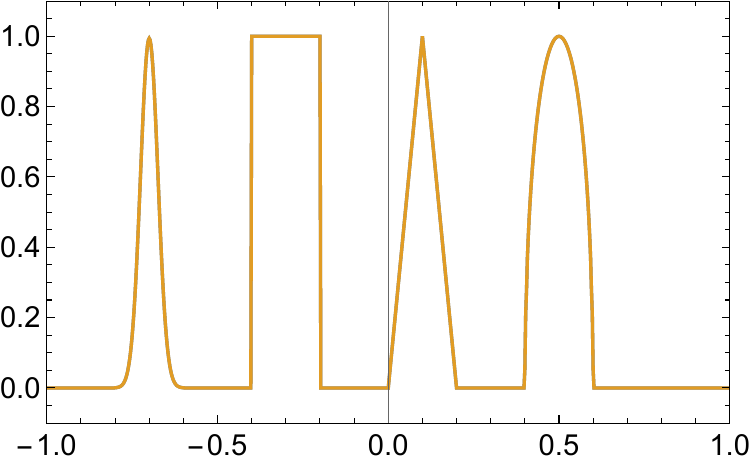}}
  \caption{
  Nonsmooth initial condition $q^0(x)$ that consists of 4 specific shapes. }
  \label{fig:init_shu}
\end{figure}

Next, we conduct a similar experiment with a non-smooth solution to evaluate the ability of the high-resolution scheme to suppress unphysical oscillations in numerical solutions. We consider a standard experiment with the initial condition $q^0(x)$ given by four different shapes as shown in Figure \ref{fig:init_shu}. 

\begin{figure}[!ht]
  \centering 
    \subfloat[3rd order accurate method]{\includegraphics[width=0.49\textwidth]{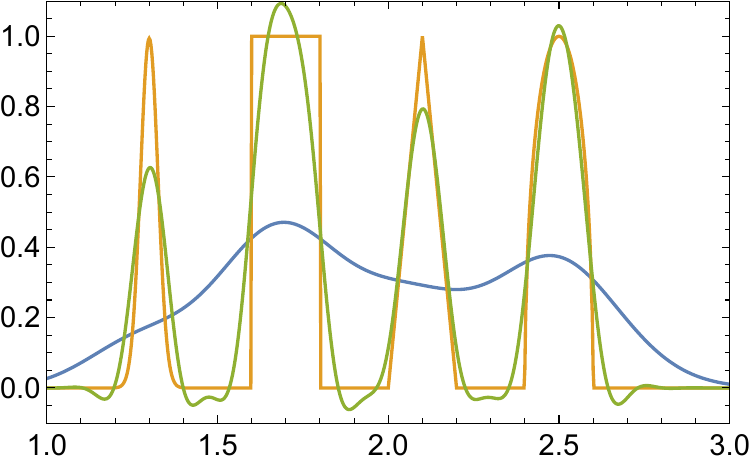}\label{fig:5a}}
    \subfloat[High resolution method]{\includegraphics[width=0.49\textwidth]{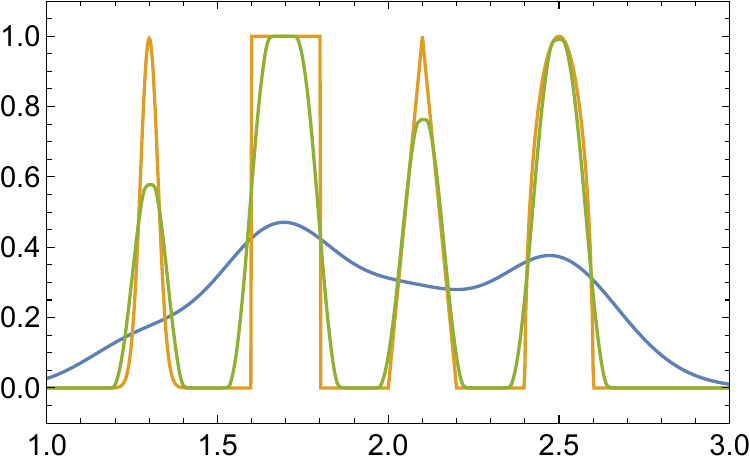}\label{fig:6a}}
  \caption{Comparison of the third-order accurate method (the green curve in left picture) and the high resolution method (the green curve in the righ picture) with the exact solution (the orange curve) and the first-order method (the blue curve) at $T=2$. The discretization step is $h=0.0025$ and the Courant number is $5$. }
  \label{fig-shu1}
\end{figure}

Each particular segment of the function $q^0(x)$ was computed separately over an interval $(a_i,3)$ with $a_i$ being $a_1=0.6$, $a_2=0.2$, $a_3=-0.2$ and $a_4=-0.6$, respectively, using the uniform steps $h=0.0025$ and $\tau=0.0125$, resulting in the Courant number of 5. The points $a_i$ coincide with the right "edge" of each shape in $q^0(x)$. Consequently, each shape enters the computational domain $(a_i,3)$ in the very first time step, and it travels the distance $2$. This approach was chosen to demonstrate the performance of the numerical methods presented for each segment in the initial profile in equal measure. The results can be seen in Figure \ref{fig-shu1}. It is apparent that the high-resolution method performs well, providing a solution with no non-physical oscillations. One additional corrector step was used in the computations only for a few values of $Q_i^n$ when the left inequality in \eqref{inequalities} was violated. In that case, no undershooting or overshooting occurred for the results even at the level of rounding errors. 

\begin{figure}[!ht]
  \centering 
    \subfloat[High-resolution method for $I=400$]{\includegraphics[width=0.49\textwidth]{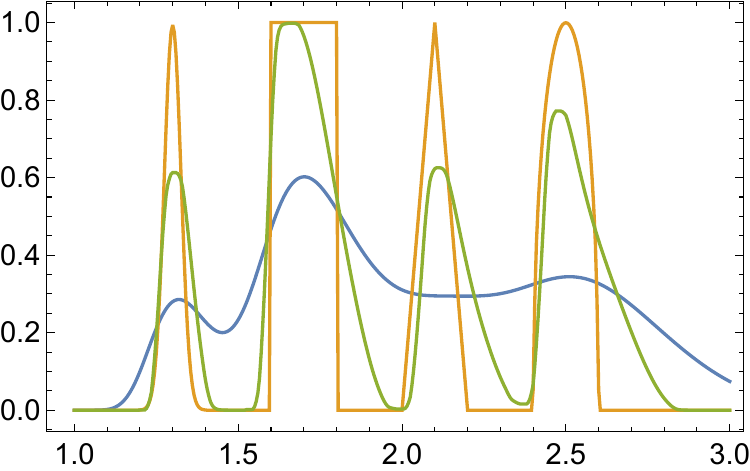}\label{fig:5}}
    \subfloat[High resolution method for $I=800$]{\includegraphics[width=0.49\textwidth]{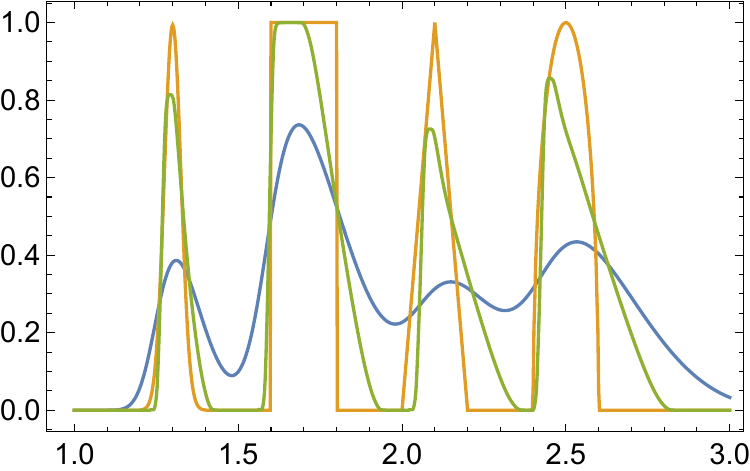}\label{fig:6}}
  \caption{Comparison of numerical solutions for nonlinear problem \eqref{nonlinear-isotherm} obtained with the high resolution method (the green curves) and the first-order method (the blue curves) at $T=2$ with the exact solution of the linear problem (the orange curves). The maximal Courant number is $\frac{10}{9}4$. }
  \label{fig-shu2}
\end{figure}

Finally, we present a nonlinear extension of this example by solving the equation,
\begin{equation}
    \label{nonlinear-isotherm}
    \partial_t \left(\frac{9}{10}q + \frac{1}{10}q^2 \right) + \partial_x q = 0 \,, \quad q(1,t)=q^0(x-t) \,,
\end{equation}
for the function $q^0$ plotted in Figure \ref{fig:init_shu}. The example is computed with the nonlinear HR version of the method as presented in Section \ref{sec-nonlinear}. The purpose is to confirm, for this example, the opinion that the HR scheme, well designed for the linear advection, can be successfully and quite straightforwardly extended to nonlinear cases. Note that $\theta'(q)=9/10+2/10 q$, so the retardation is highest for $q=1$ and lowest for $q=0$. One can see in Figure \ref{fig-shu2} for two different grids with the maximal Courant number $C^{\max}=\frac{10}{9}4$ that no under- and overshooting occurs in numerical solutions that can be confirmed also numerically even at the level of rounding errors. The same strategy of one additional corrector step as described for the linear case was used for $C^{\min}=\frac{10}{11} 4$.

One can clearly see in Figure \ref{fig-shu2} the typical features of nonlinear scalar hyperbolic equations in the form of rarefaction waves and shock waves \cite{frolkovicSemianalyticalSolutionsContaminant2006}. 

\subsection{Benchmark on a triangular network}
\label{sec:triangle}

In this section, we compare our method with the results on a network as proposed in \cite{eimerImplicitFiniteVolume2022, barsukowImplicitActiveFlux2024}. The purpose is to assess the accuracy and stability of the proposed schemes for medium Courant numbers for a nontrivial form of the network. 

\begin{figure}[!ht]
  \centering
    \includegraphics[width=0.55\textwidth]{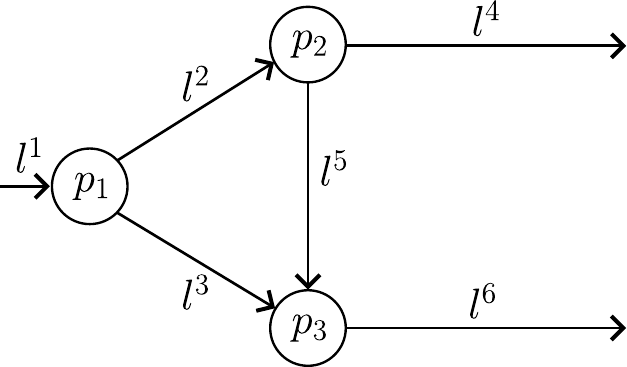}
    \label{fig:BBnetwork}
  \caption{Triangular network proposed by Barsukow and Borche \cite{barsukowImplicitActiveFlux2024}}
\label{fig-barsukow}
\end{figure}

The network consists of 3 inner vertices $p_1,\ p_2,\ p_3$ and 6 edges $l^1,\ldots,l^6$, as shown in Figure \ref{fig-barsukow}. Note that in order to be consistent with the description in \cite{eimerImplicitFiniteVolume2022,barsukowImplicitActiveFlux2024}, we do not denote the inflow vertex and two outflow vertices in this figure. The lengths of the edges are given by $L^1=5$, $L^2=L^3=L^5=20$, $L^4=L^6=30$ and the velocities are constant per edge with $v^1=v^3=v^4=v^6=1$,  $v^2=2$ and $v^5=\frac{40}{23}$.

Considering the boundary conditions, the flux is split in the vertex $p_1$,
\begin{eqnarray}
    \label{alpha divide2}
    q^2(t,0) = \alpha_1 q^{1}(t,L^1)\,, \quad
    q^{3}(t,0) = (1-\alpha_1) q^{1}(t,L^1) \,,  
\end{eqnarray}
and $p_2$,
\begin{eqnarray}
    \label{alpha couple}
    q^{4}(t,0) = \alpha_2 q^{2}(t,L^2)\,, \quad
    q^{5}(t,0) = (1-\alpha_2) q^{2}(t,L^2) \,.  
\end{eqnarray}
The constants take the values $\alpha_1=\frac{3}{4}$ and $\alpha_2=\frac{2}{3}$.
For the vertex $p_3$ the coupling is simply defined by
\begin{eqnarray}
    \label{alpha divide}
    q^{6}(t,0) = q^{3}(t,L^3)+q^{5}(t,L^5).  
\end{eqnarray}
The inflow boundary condition is given by
\begin{eqnarray}
    \label{BB_BC}
    q_{e_1}(t,0) = \sin\Big( \frac{2\pi}{3}t\Big) \,.
\end{eqnarray}
To be consistent with \cite{barsukowImplicitActiveFlux2024}, we define also nonzero initial conditions only
for the first edge,
\begin{eqnarray}
    \label{BB_IC}
    q^{1}(t,0) =  e^{-4\big(x-2.5\big)^2}.
\end{eqnarray}

\begin{figure}[H]
  \centering
  \subfloat[t = 0]{\includegraphics[width=0.44\textwidth]{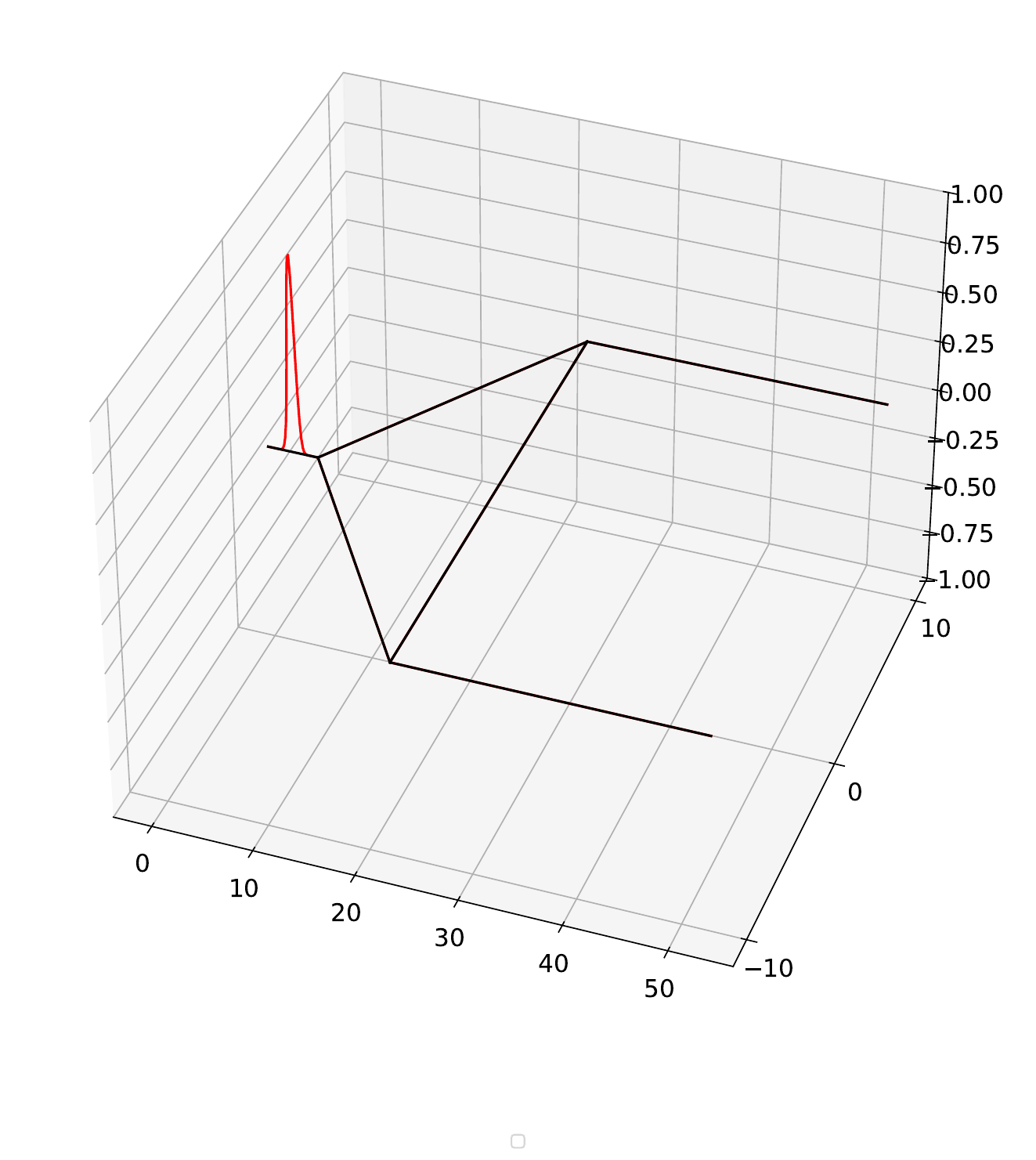}\label{fig:BB_noENO3d_t0}}
  \subfloat[t = 15]{\includegraphics[width=0.44\textwidth]{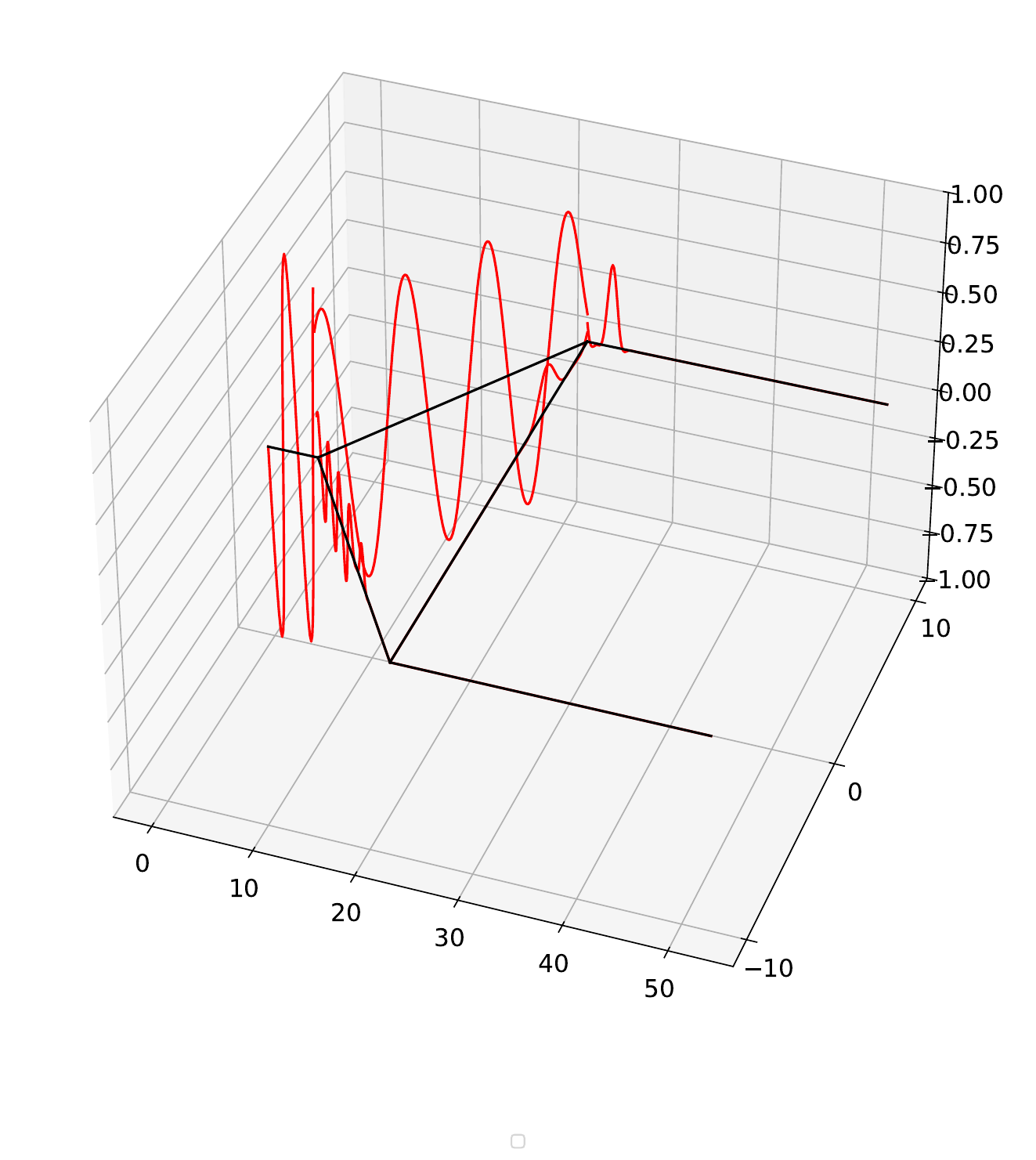}\label{fig:BB_noENO3d_t15}}

  \subfloat[t = 30]{\includegraphics[width=0.44\textwidth]{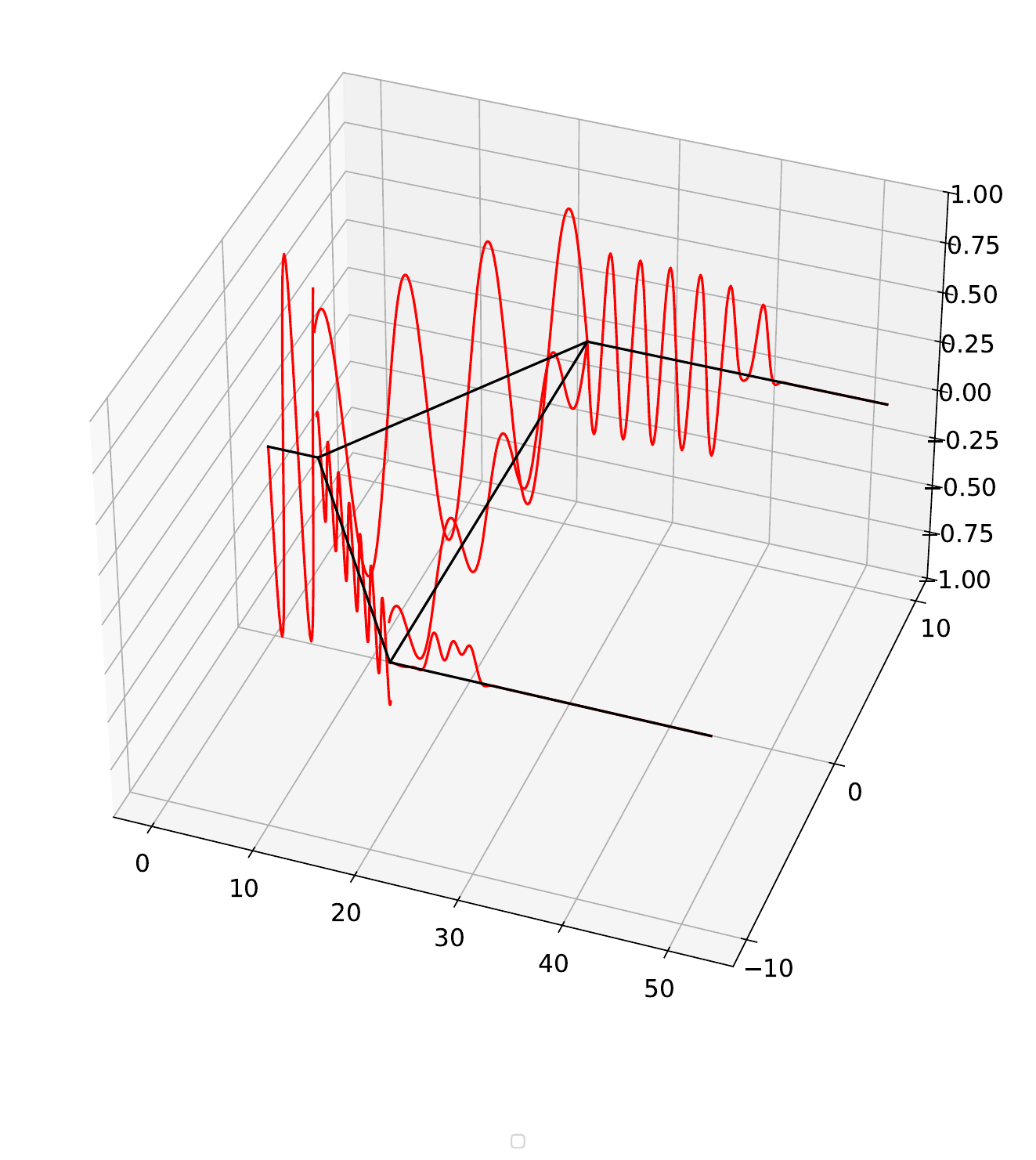}\label{fig:BB_noENO3d_t30}}
  \subfloat[t = 70]{\includegraphics[width=0.44\textwidth]{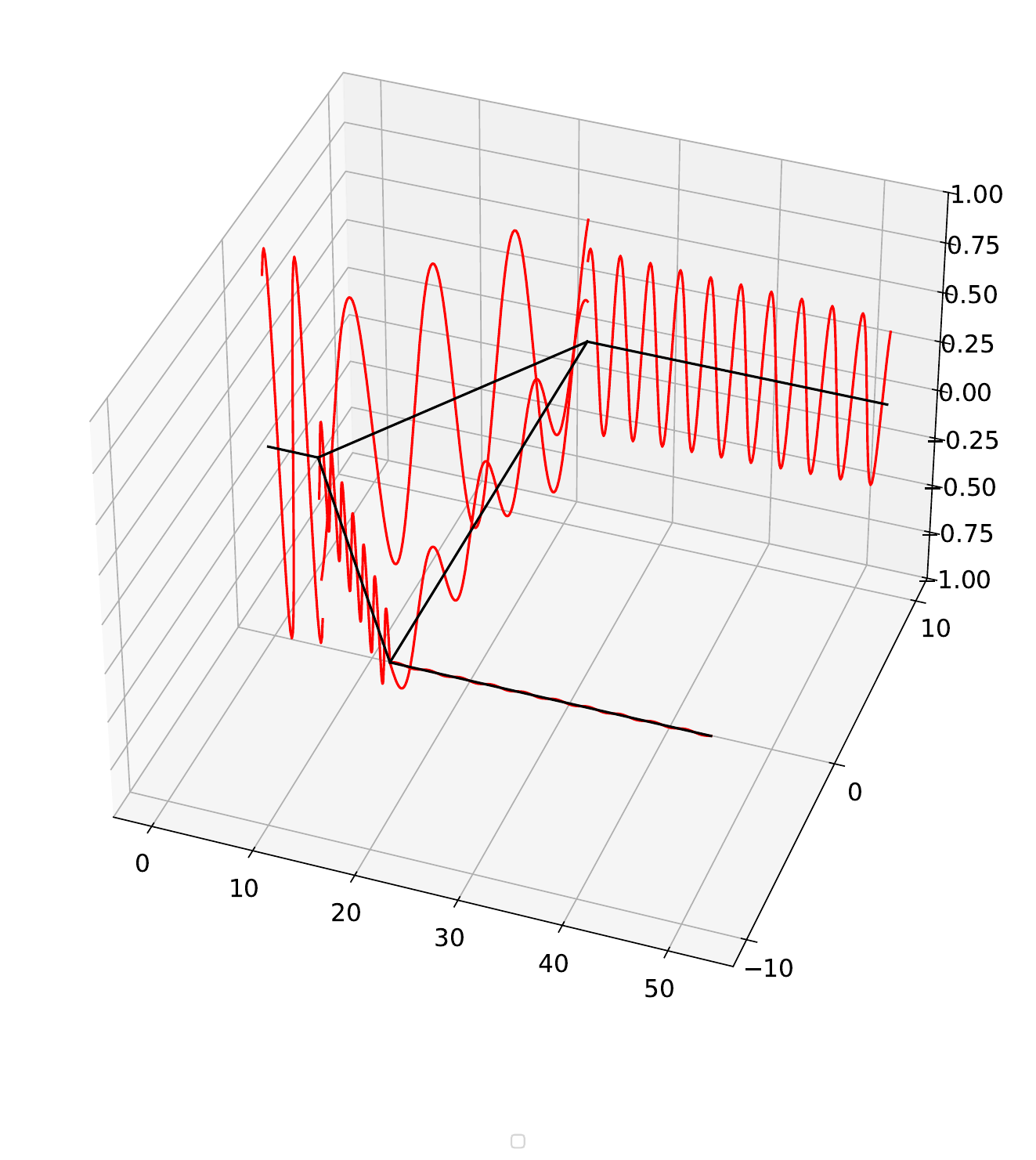}\label{fig:BB_noENO3d_t70}}

  \caption{Benchmark test on a triangular network. The numerical solutions are obtained with the third-order accurate method.}
  \label{fig:BB_noENO3d}
\end{figure}

First, we present the initial condition and numerical solutions computed with the third order accurate scheme at three different times in Figure \ref{fig:BB_noENO3d}. The splitting of the initial Gaussian and the sine wave can be observed in Figure \ref{fig:BB_noENO3d_t15} at $t=15$. The experiment is designed in such a way that the exact solution for large times is equal to zero at the edge $l^6$, see Figure \ref{fig:BB_noENO3d_t70}. This is caused by the interference of three incoming waves, two originating from the initial Gaussian wave and one originating from the sine wave at the vertex $p_2$, as observed in Figure \ref{fig:BB_noENO3d_t30}.

To achieve a comparable accuracy in this benchmark experiment as in \cite{barsukowImplicitActiveFlux2024}, the results were obtained using the following discretization parameters. We choose the same space discretization step as in \cite{barsukowImplicitActiveFlux2024}, namely $h=1/8$. In \cite{barsukowImplicitActiveFlux2024} the time step is chosen $\tau=5/4$ that corresponds to the maximal Courant number equal to $C_{\max}=5$. The method used in \cite{barsukowImplicitActiveFlux2024} is the fourth-order accurate Active Flux finite-volume method with unknowns not only in the cells represented by nodes $x_i$, but also in the middle points $x_{i+1/2}=x_i + h/2$. Furthermore, the updates of numerical values are more complex than in our scheme as they contain, for instance, coefficients with the $5$th order polynomials of $C$. The scheme used in \cite{barsukowImplicitActiveFlux2024} is stable only for $C>1.1$. 

We do not aim to compare the efficiency of these two schemes, but we can clearly compensate for the advantages of the scheme in \cite{barsukowImplicitActiveFlux2024} by using finer discretization steps with our third-order scheme.
As precision is the main objective of this experiment, we choose to compute the example with the time step $\tau=5/16$ that results in constant Courant numbers $C^e$ per each edge having values $C^1=C^3=C^4=C^6=0.625$, $C^2=1.25$, and $C^5\approx 1.087$. For such choices of $h$ and $\tau$, the results in Figure \ref{fig:BB_noENO3d} compare well with those presented in \cite{barsukowImplicitActiveFlux2024}. 

The accuracy of the numerical results is critical for the approximations of numerical solution on the sixth edge. As described above, due to the interference of the incoming waves, the exact solution should be zero; see Figure \ref{fig:BB_noENO3d_t70}, except for the three isolated waves, see Figure \ref{fig:BB_noENO3d_t30}. Therefore, analogously to \cite{barsukowImplicitActiveFlux2024}, we provide a numerical study of the convergence for the approximation of $q^6$ for $t=35$ and $t=70$.

\begin{figure}[H]
  \centering
    \subfloat[Results on a coarser mesh]{\includegraphics[width=0.49\textwidth]{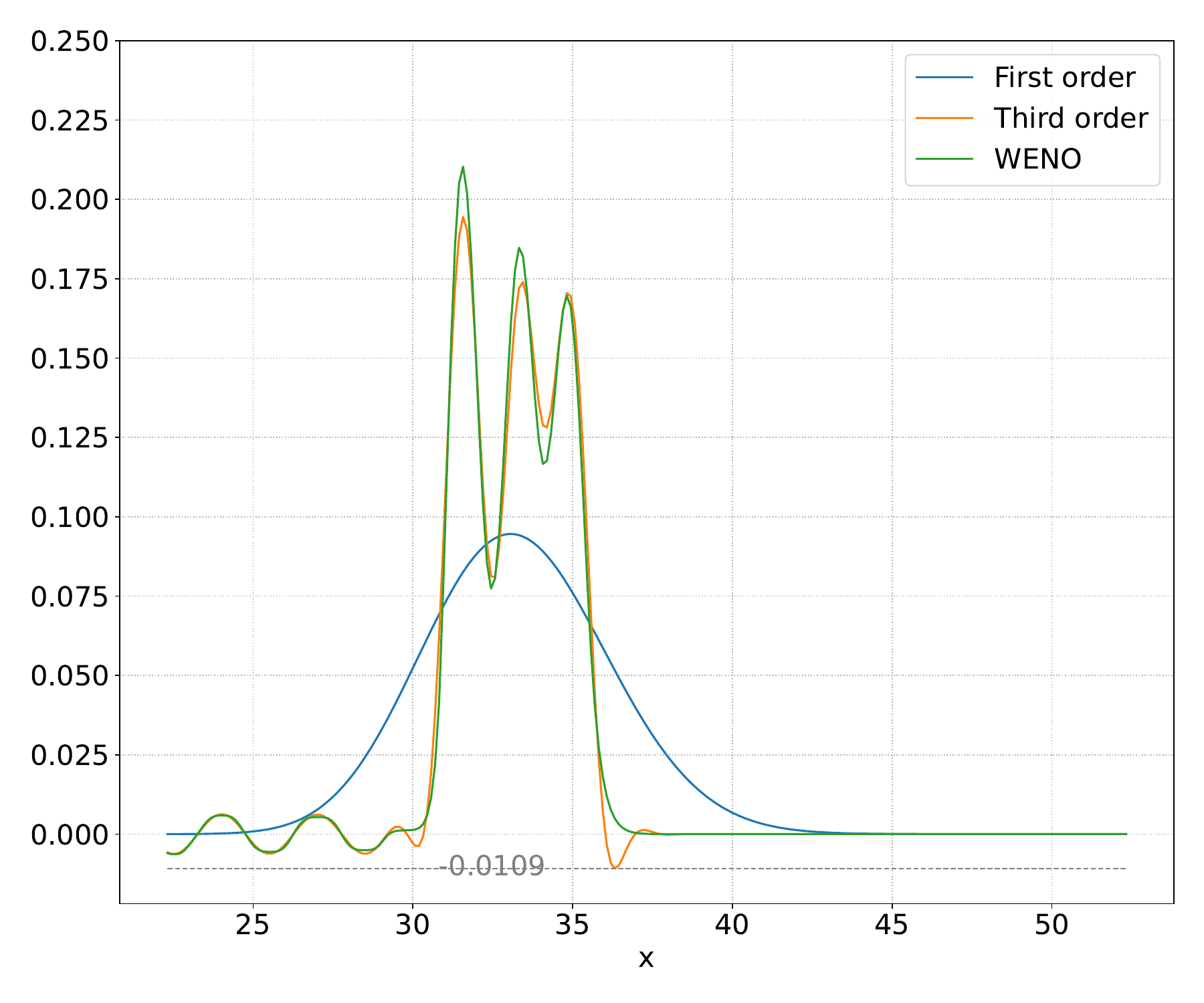}\label{fig:bbnoeno35}}
    \subfloat[Results on a finer mesh]{\includegraphics[width=0.49\textwidth]{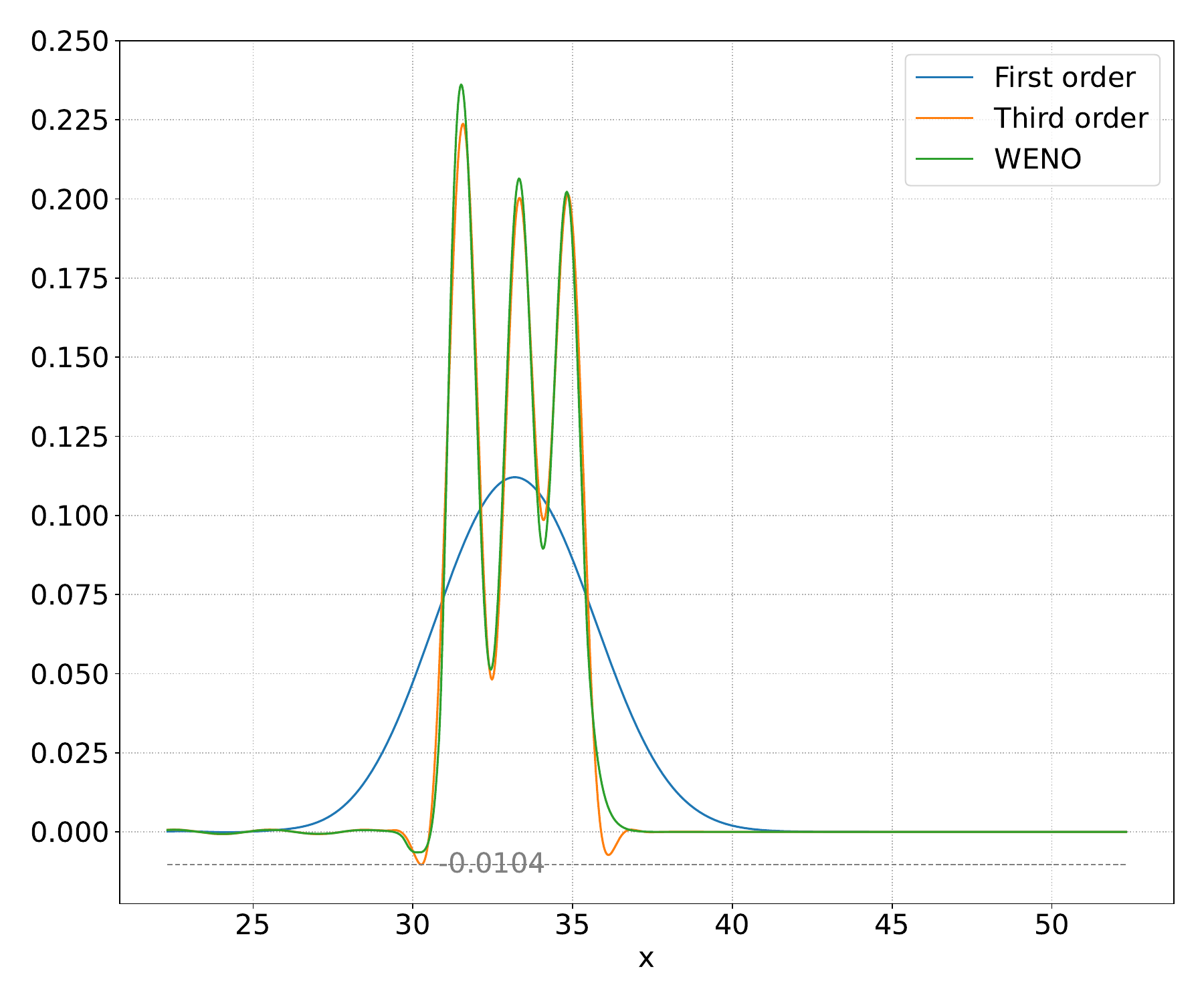}\label{fig:bbnoeno35finer}}
  \caption{Comparison of the third-order scheme and the WENO for the numerical approximation of $q^6(x,35)$ using the same time step $\tau =5/16$ and different Courant numbers. In Figure \ref{fig:bbnoeno35} the values are $C^1=C^3=C^4=C^6=0.625$, $C^2=1.25$ and $C^5\approx 1.087$ and in Figure \ref{fig:bbnoeno35finer} the values are $C^1=C^2=C^3=C^4=C^6=2.5$, and $C^5\approx 2.174$.}
  \label{fig:BB_compare_35}
\end{figure}

Concerning the results for $t=35$ in Figure \ref{fig:BB_compare_35}, we present the numerical solutions obtained with the same settings as in Figure \ref{fig:BB_noENO3d}. In contrast to the first-order scheme, an approximation of all three waves can be observed in Figure \ref{fig:bbnoeno35}. However, the accuracy is reduced because of a smoothing of local extremes and because of unphysical oscillations represented by negative values. Therefore, we also present the numerical solutions computed with the same settings but using the WENO scheme from Remark \ref{rem-weno} to reduce the negative oscillations. In Figure \ref{fig:bbnoeno35} the benefits of the WENO scheme can be clearly seen.

To further improve the accuracy of the results and, at the same time, to illustrate a stable behavior of numerical schemes for larger Courant numbers, we keep the time step $\tau=5/16$ and refine each edge so that the constant Courant numbers at each edge now take the values $C^1=C^2=C^3=C^4=C^6=2.5$, and $C^5\approx 2.174$. In Figure \ref{fig:bbnoeno35finer}, we compare the results for the third-order scheme and the WENO scheme. Both numerical solutions clearly contain less smoothing of the local extrema, and the WENO scheme reduces the negative oscillations. Note that in Figure \ref{fig:bbnoeno35}, to avoid smoothing of local extremes, the WENO scheme is applied only for the predicted values $|\widetilde{Q}_i^n|<0.1$ and in Figure \ref{fig:bbnoeno35finer} for $|\widetilde{Q}_i^n|<0.05$. 

\begin{figure}[H]
  \centering
    \subfloat[Results on the coarser mesh]{\includegraphics[width=0.49\textwidth]{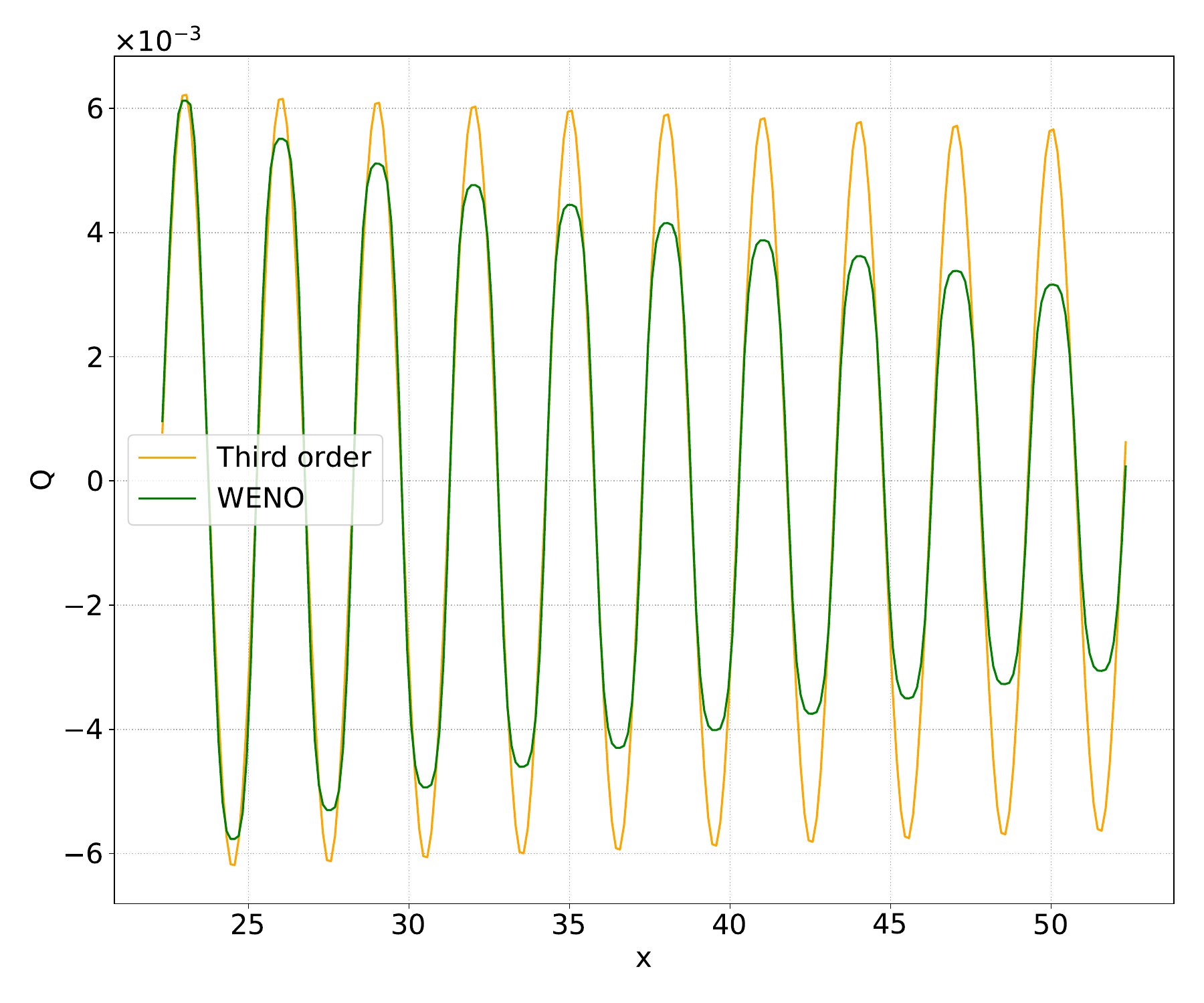}\label{fig:bbnoeno70}}
    \subfloat[Results on the finer mesh]{\includegraphics[width=0.49\textwidth]{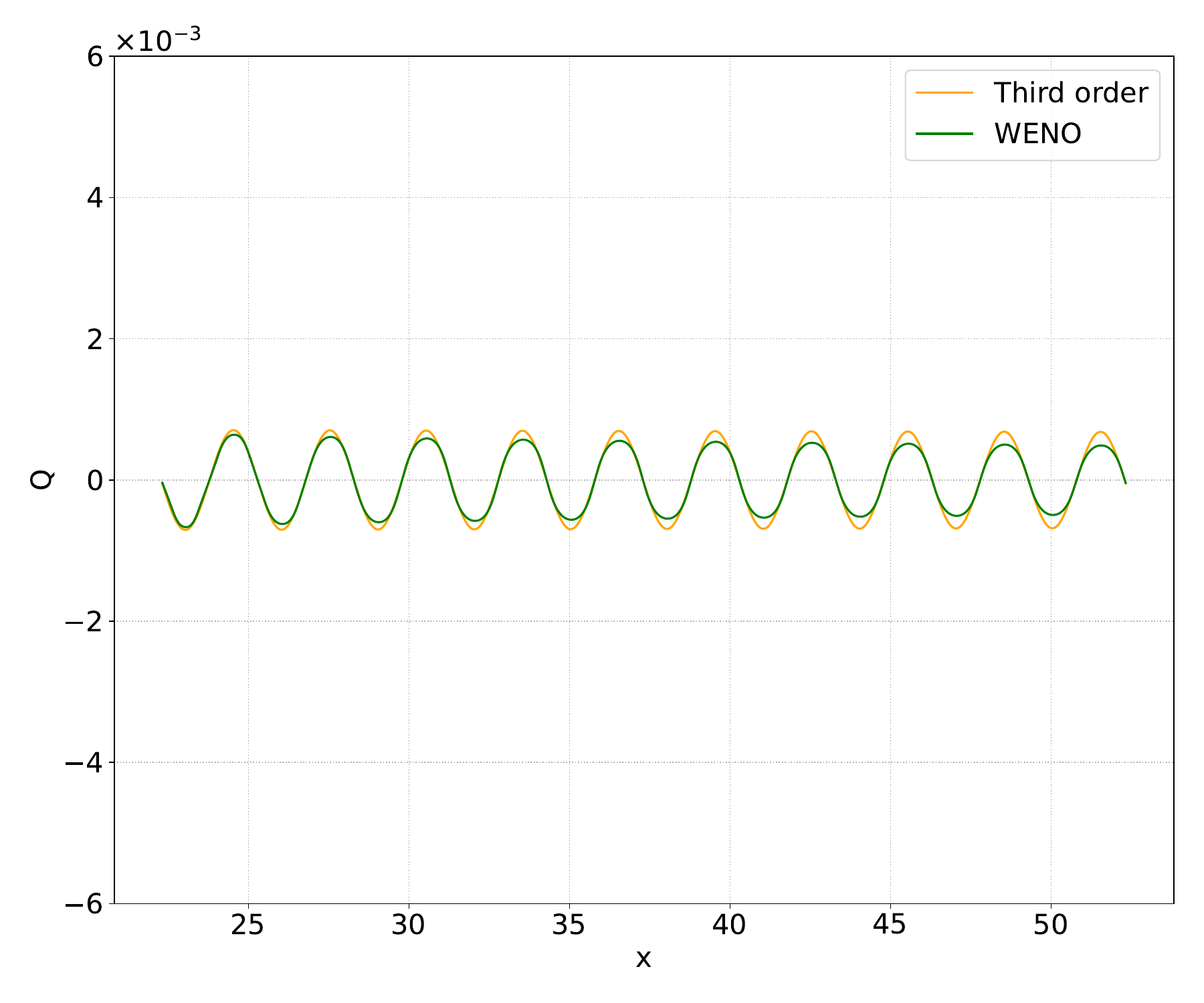}\label{fig:bboriginal70}}
  \caption{Comparison of the third-order scheme and the WENO for the numerical approximation of $q^6(x,70)$ using the $\tau =5/16$ and $h=1/8$ in Figure \ref{fig:bbnoeno70} and $\tau =5/32$ and $h=1/16$ in Figure \ref{fig:bboriginal70}.}
  \label{fig:BB_compare_70}
\end{figure}

Finally, we compare the results for $t=70$ when the exact solution $q^6$ is zero. In Figure \ref{fig:BB_compare_70}, we see that the numerical solutions take the form of a sine wave with a very small amplitude that is significantly reduced for the finer grid in Figure \ref{fig:bboriginal70} with respect to the results on the coarser mesh in Figure \ref{fig:bbnoeno70}. 

We can summarize that, for this example, we obtain satisfactory results when compared with the more involved numerical scheme used in \cite{barsukowImplicitActiveFlux2024}. The advantages of our schemes are unconditional stability, the bounds preservation and simplicity.


\subsection{Numerical experiments on a realistic network}
\label{sec-revuca}

\begin{figure}[H]
  \centering
   {\includegraphics[width=0.45\textwidth]{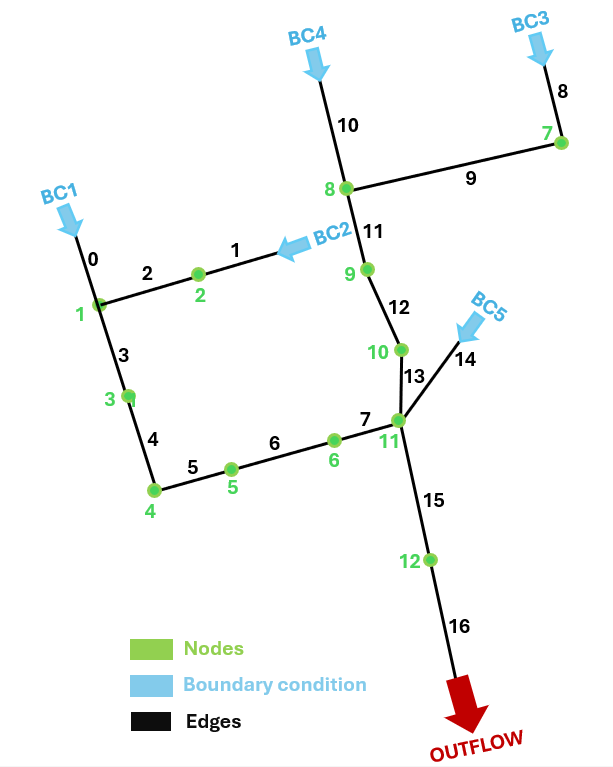}}
  \caption{Two dimensional representation of the part of the sewer network in Revúca. }
  \label{fig:sewer2D}
\end{figure}

In this section, we apply our numerical schemes to a part of a real sewer network in the city of Revúca. We acknowledge that the data were kindly provided by Doc. Ing. Marek Sokáč, CSc. from the Institute of Hydrology of the Slovak Academy of Sciences in Bratislava. The real network is simplified by shortening some parts that do not contain any branching of pipes. The final form of the network with which we are working is plotted in Figures \ref{fig:sewer2D} and \ref{fig:sewer3D}. We note that in the current version we derive a nondimensional model based on this network to assess the behavior of our numerical method for some possible scenarios of pollutant transport \cite{sokacImpactSedimentLayer2021} in complex networks.  

The network consists of five input vertices that are indicated in Figure \ref{fig:sewer2D} from \texttt{BC1} up to \texttt{BC5} and $12$ inner vertices numbered from \texttt{1} to \texttt{12}. In addition, there is one output vertex denoted by \texttt{OUTFLOW}. All vertices in the provided data are specified not only by their coordinates in the plane but also by their depth. The vertices are connected by 17 edges that are numbered from \texttt{0} to \texttt{16}. The numbering of vertices and edges follows the rules described in Section \ref{sec-mathematical-model}.

Furthermore, each edge is characterized by a constant diameter of the cross-sectional area that is variable throughout the network. In our non-dimensional settings, the values of the diameter vary from $0.3$ to $1$. The lengths of the edges are normalized, so their non-dimensional lengths vary from approximately $0.24$ to $1$.
The flow direction clearly follows the edges from the input vertices to the single output vertex. This behavior is supported by the data, as the depth of the vertices decreases along the flow direction. In addition, we computed the elevation angles of each edge based on the depth coordinates of the vertices.

There are three coupling vertices, two edges cross at the vertex \texttt{1} and \texttt{8}, and the vertex \texttt{11} is a cross-section point of three edges. The coupling boundary conditions for these three vertices are given by \eqref{incomming-concentration}, for all the other vertices, the simple continuity condition \eqref{simple-coupling} is used.

In our experiment, we suppose constant velocities per each edge to represent a stationary flow in the network. The values are obtained by prescribing constant nondimensional water discharge (volume flow) for each inflow vertex. Due to the varying diameters and inclinations of the pipes, the velocities differ across each edge, as they are computed using the continuity equation.

\begin{figure}[H]    
  \centering
    \includegraphics[width=0.6\textwidth]{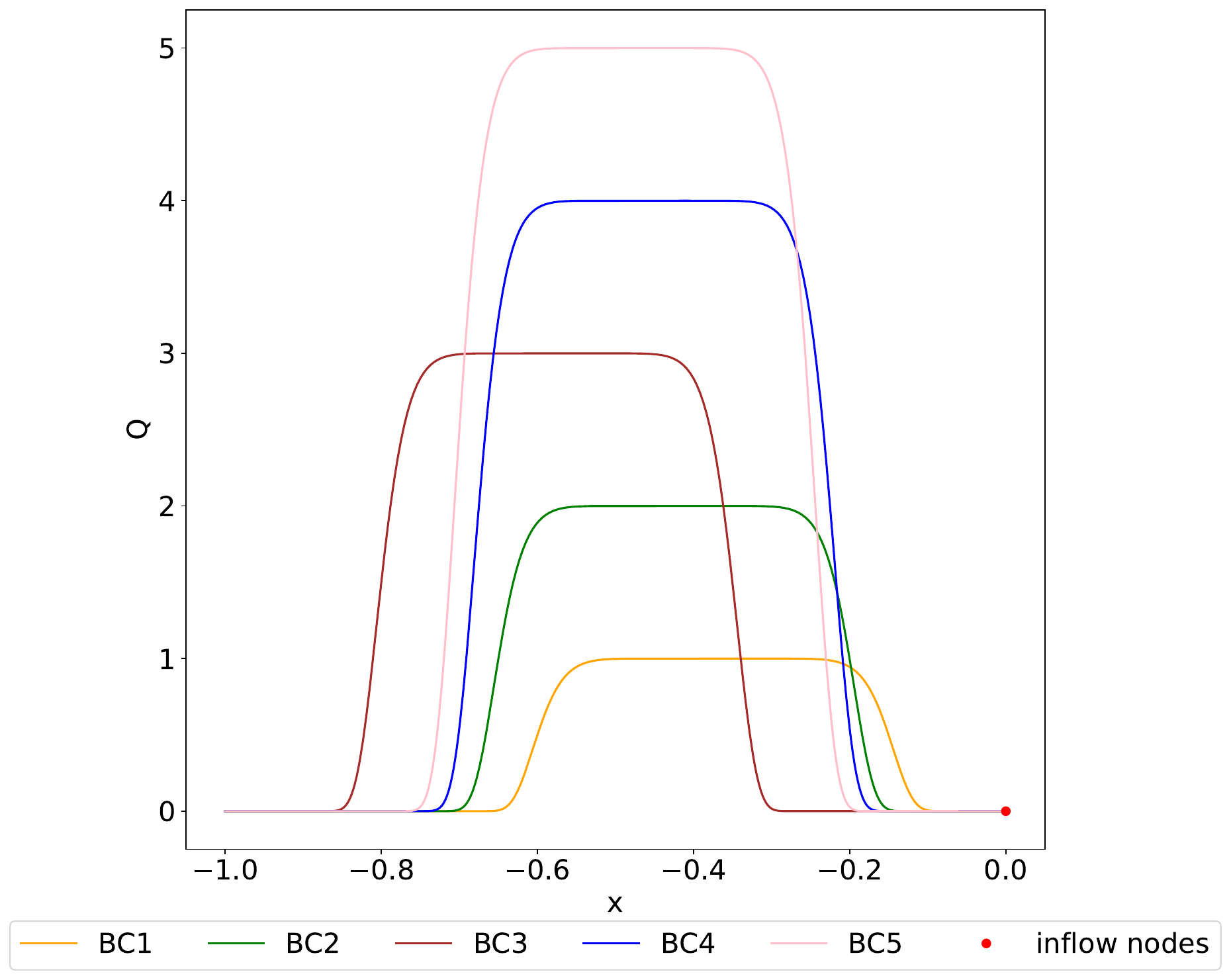}
  \caption{Five specific boundary inflow functions \eqref{boundary_function} at time $t=0$.}\label{fig:boundarygraph}
\end{figure}

The boundary condition for the inflow vertices are prescribed by
\begin{eqnarray}
    \label{boundary_function}
    q_m(t) = q^{in}_m(x,t) = c_m e^{-2(3(x-0.3t)-x_m)^{10}} \,, c_m>0 \,, x_m<0 \,,
\end{eqnarray}
for $m$ from \texttt{BC1} to \texttt{BC5}. The parameters $c_m$ and $x_m$ were chosen so that each function $q^{in}_m$ would represent an "impulse" (a profile consisting of a sharp increase from 0 to the maximal value, followed by an almost constant plateau and a sharp decrease to zero) with a different maximal value and a slightly different initial position. The maximum height of the wave corresponds to the index number of the vertex $m$. The particular functions $q^{in}_m$ in $t=0$ are plotted in Figure \ref{fig:boundarygraph}.  

For our numerical experiment in this network, we set the final simulation time $T = 2$. The number of time steps is $N = 384$. The number of spatial steps $I^e$ per each edge $l^e$ is 128 or 256 to take into account the significantly different lengths of the edges in the network.
The Courant numbers $C^e$ for each edge are shown in Table \ref{tab:cournatnumbs_half}. 

\begin{table}[!ht]
\footnotesize
    \centering
    \begin{tabular}{c|c|c|c|c|c|c|c|c|c}
    \hline
        $e$ & 0 & 1 & 2 & 3 & 4 & 5 & 6 & 7 & 8 \\ \hline
            \hline 
        $I^e$ & 128 & 128 & 128 & 128 & 128 & 128 & 128 & 128 & 128 \\ \hline
        $C^e$ & 2.61 & 1.16 & 1.21 & 3.70 & 3.58 & 2.37 & 1.69 & 2.56 & 8.015 \\ \hline
        \multicolumn{10}{c}{ }  \\
        \hline
        $e$ & 9 & 10 & 11 & 12 & 13 & 14 & 15 & 16 & ~ \\ \hline\hline
        $I^e$ & 256 & 128 & 128 & 128 & 128 & 128 & 256 & 256 & ~ \\ \hline
        $C^e$ & 0.88 & 2.12 & 5.37 & 5.18 & 6.01 & 2.42 & 3.52 & 2.15 & ~ \\ \hline
    \end{tabular}
    \caption{Courant numbers for edges, N = 384.}
    \label{tab:cournatnumbs_half}
\end{table}

\begin{figure}[H]
  \centering
  \subfloat[Path $BC1$]{\includegraphics[width=0.43\textwidth]{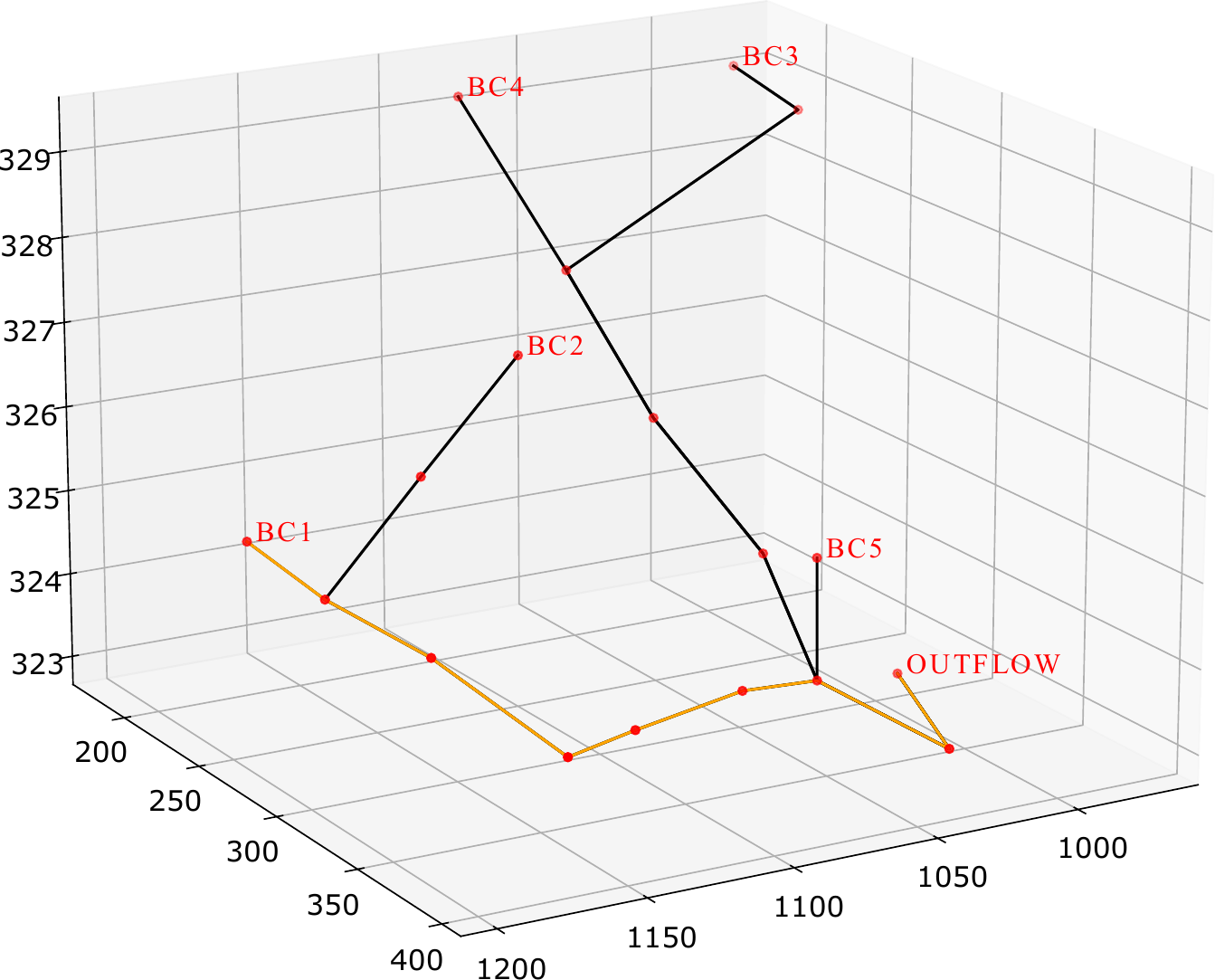}\label{fig:f_way1}}
  \subfloat[Path $BC3$]{\includegraphics[width=0.43\textwidth]{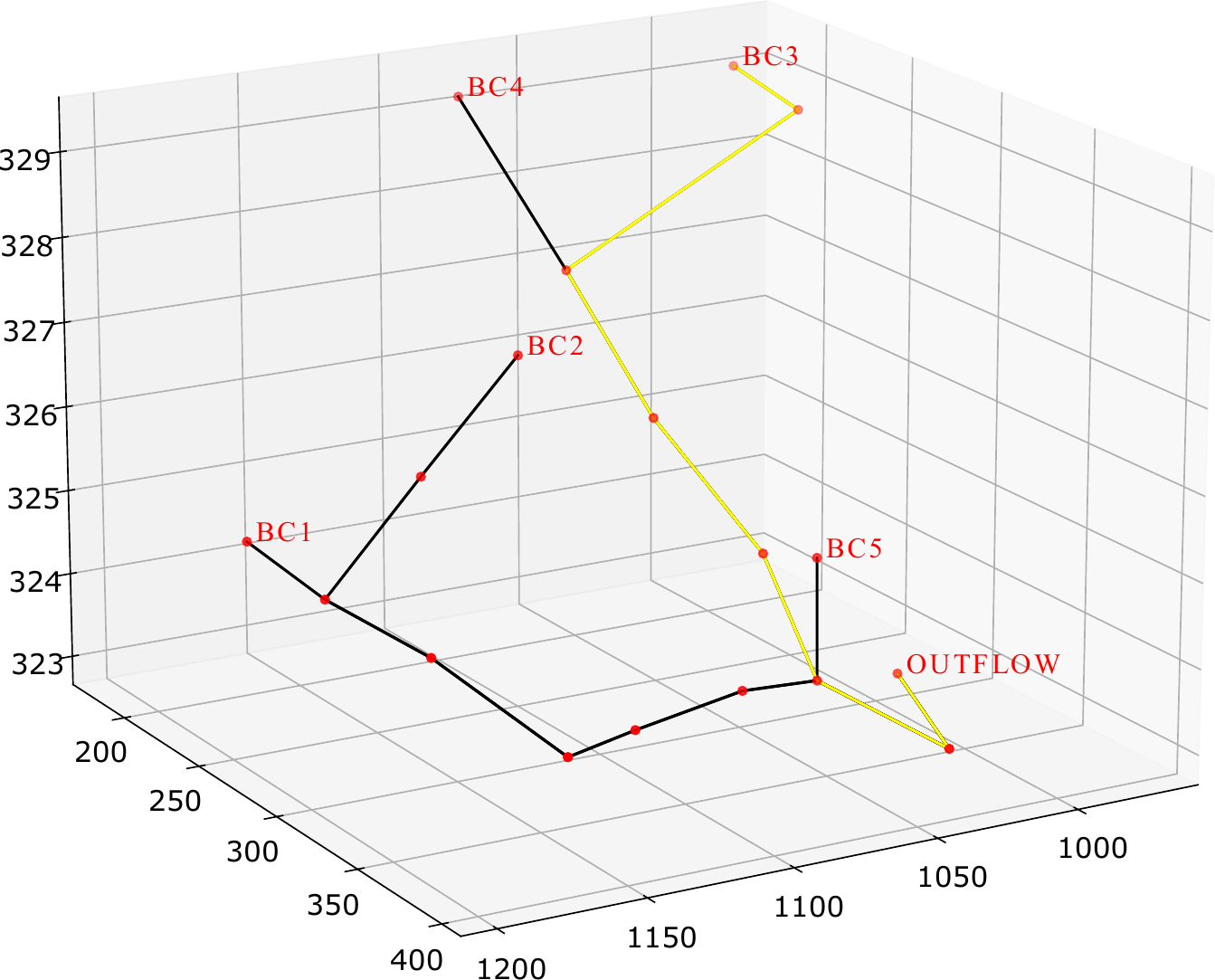}\label{fig:f_way3}}
  \caption{Two selected paths of the Revúca sewer network.}
  \label{fig:f_ways}
\end{figure}

To demonstrate the results, we mark two paths originating from the inflow vertices $BC1$ and $BC3$, respectively, as shown in Figures \ref{fig:f_ways}. The naming of the paths is based on the inflow vertex from which each path originates, with each path terminating at the \texttt{OUTFLOW} vertex. Note that the part after the coupling vertex 11 is common for both presented paths.

\begin{figure}[H]
  \centering
  \subfloat[t = 0.25]{\includegraphics[width=0.44\textwidth]{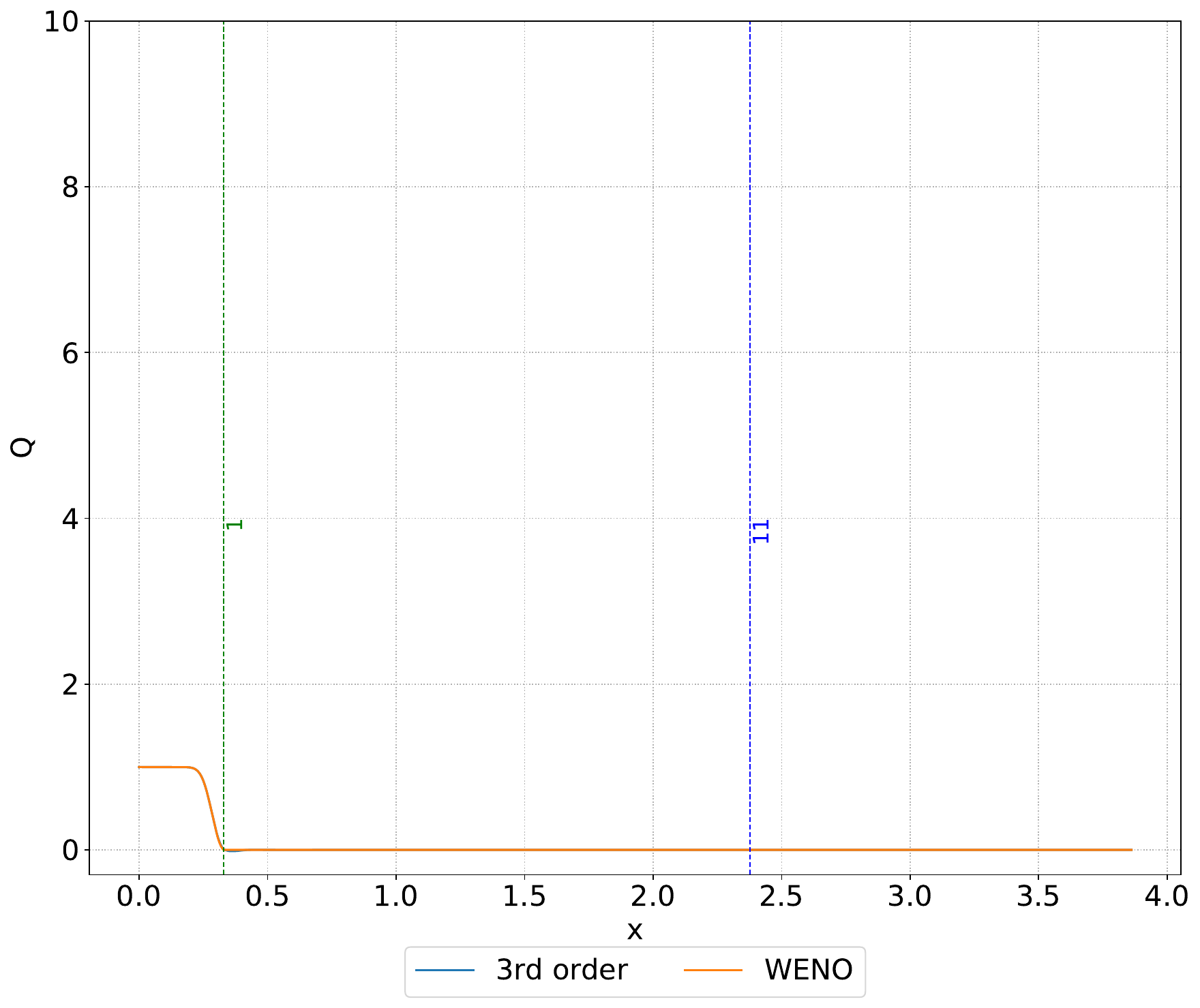}\label{fig:f_bc1_0125t_halfCournat}}
  \subfloat[t = 0.5]{\includegraphics[width=0.44\textwidth]{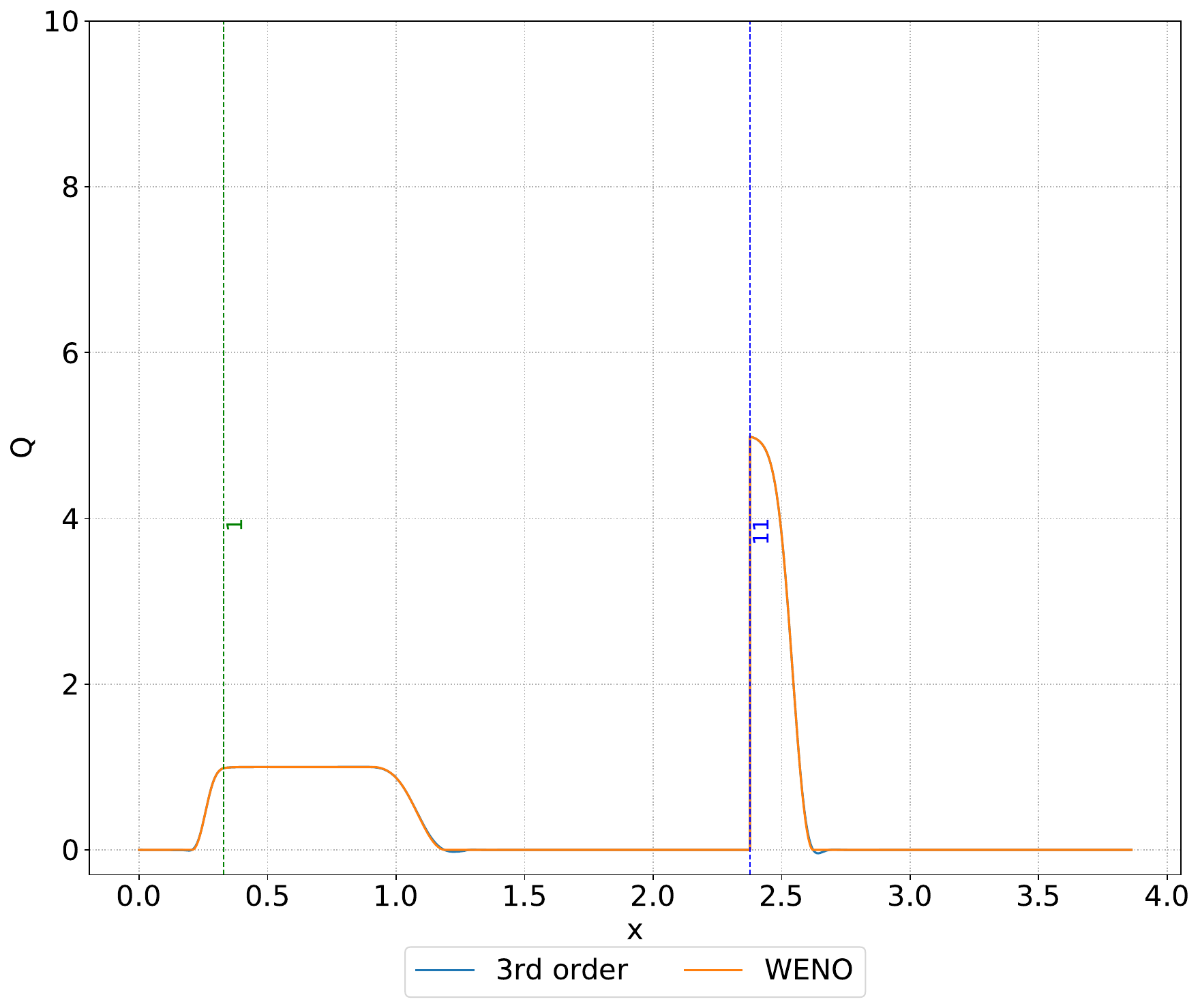}\label{fig:f_bc1_025t_halfCournat}}

  \subfloat[t = 1.0]{\includegraphics[width=0.44\textwidth]{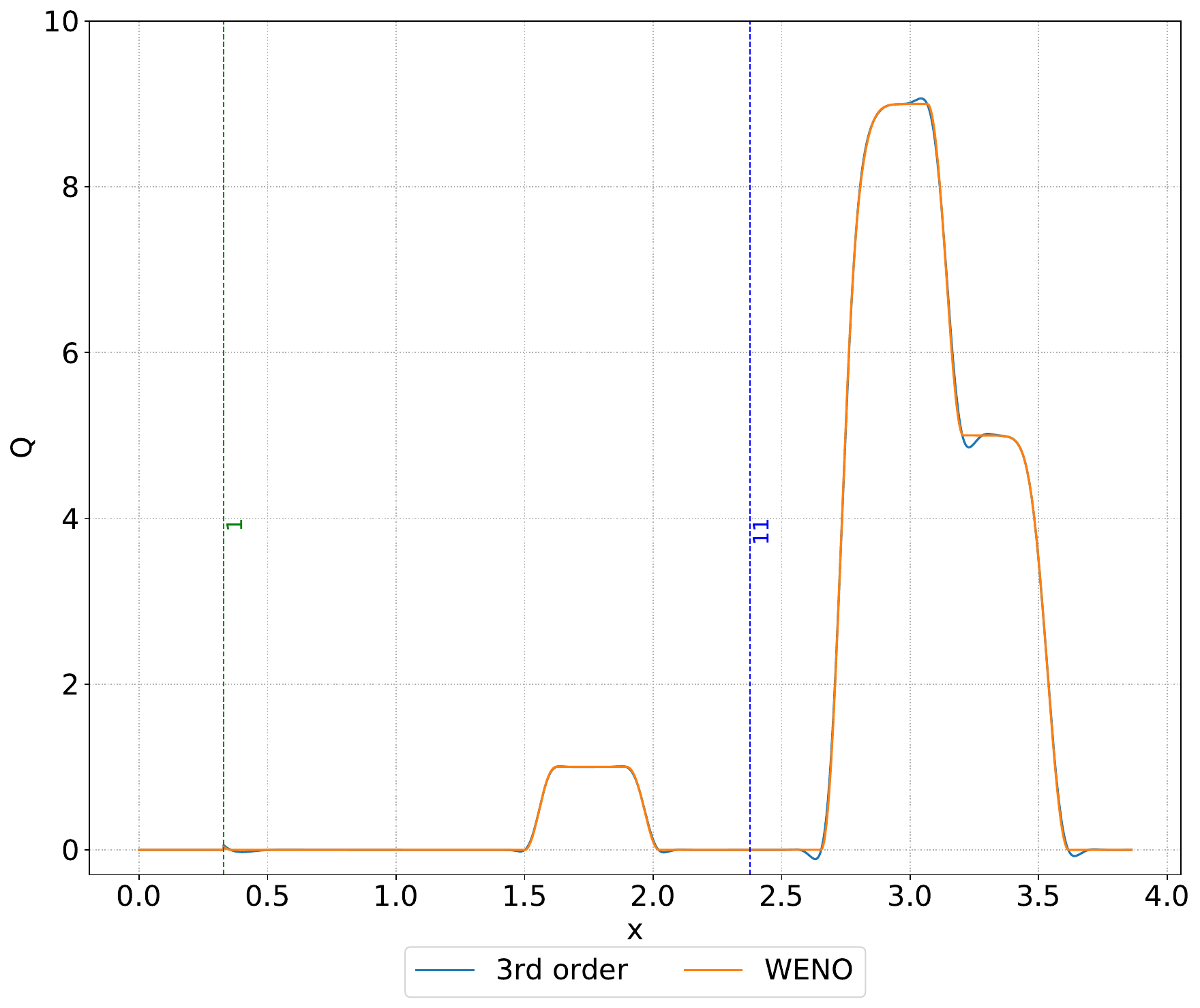}\label{fig:f_bc1_05t_halfCournat}}
  \subfloat[t = 2.0]{\includegraphics[width=0.44\textwidth]{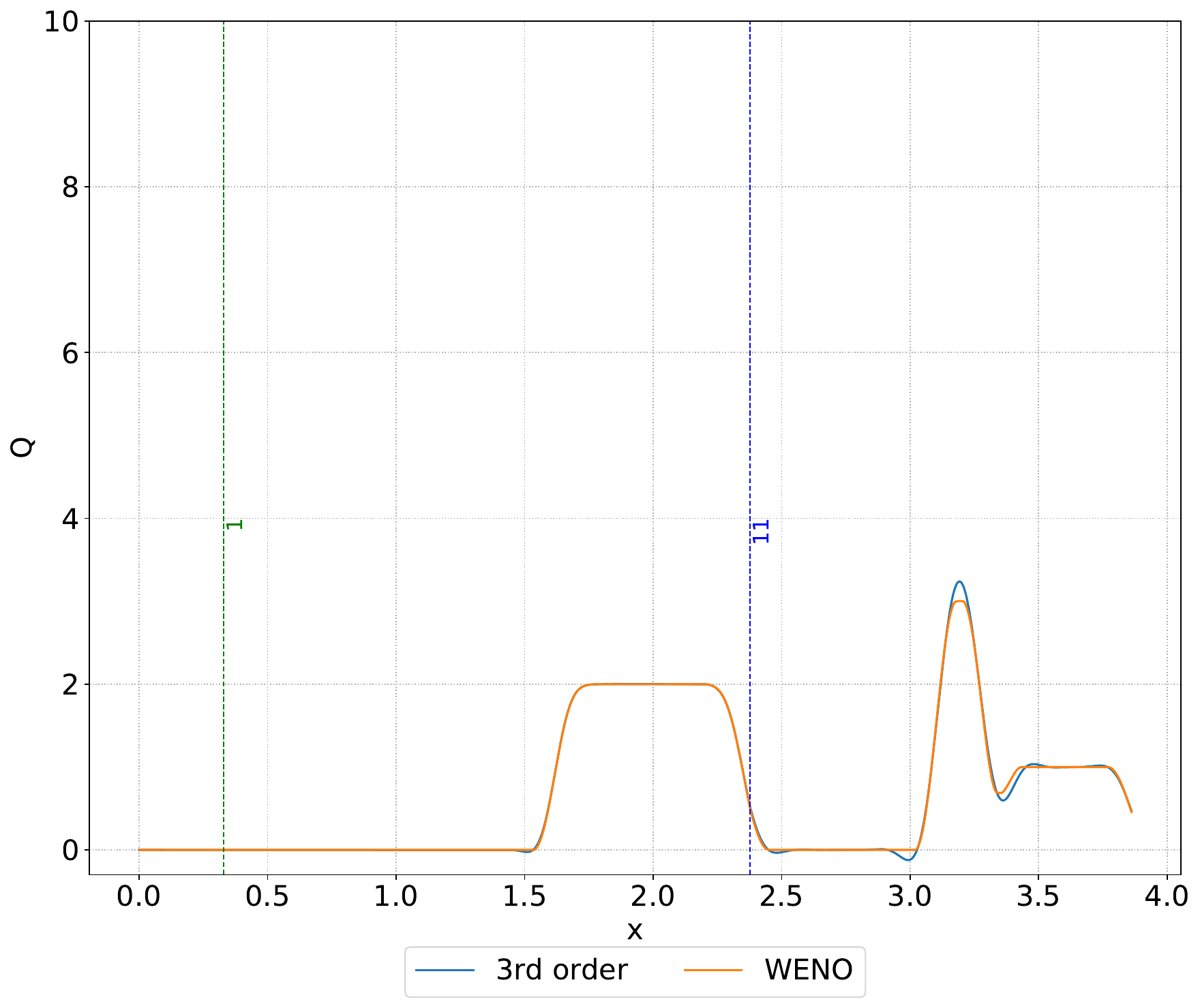}\label{fig:f_bc1_1t_halfCournat}}

  \caption{Flow through the path $BC1$ shown in 4 different moments. Courant numbers on each edge are shown in Table \ref{tab:cournatnumbs_half}.}
  \label{fig:f_way1_timeT_halfCournat}
\end{figure}

\begin{figure}[H]
  \centering
  \subfloat[t = 0.25]{\includegraphics[width=0.44\textwidth]{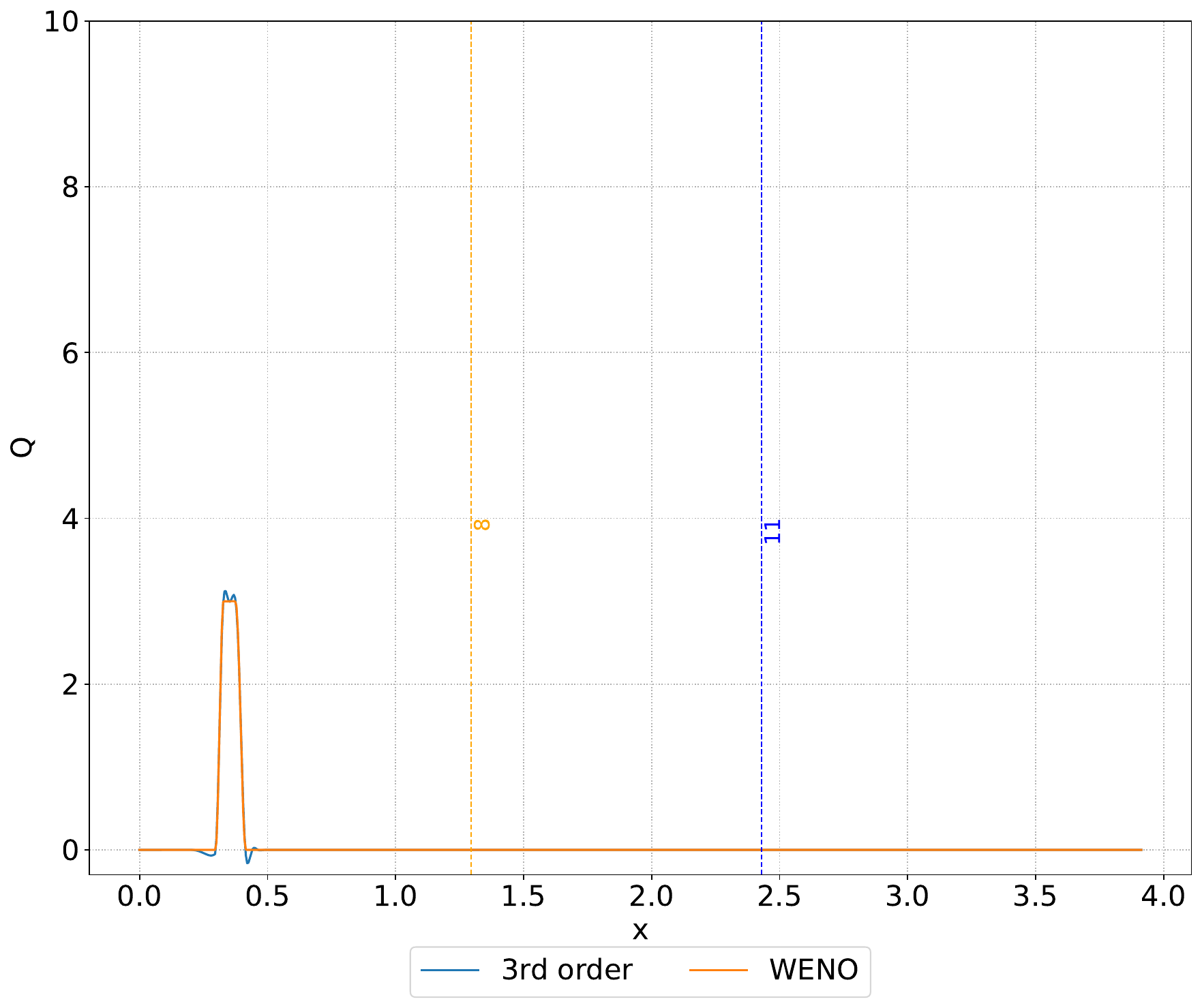}\label{fig:f_bc3_0125t_halfCournat}}
  \subfloat[t = 0.5]{\includegraphics[width=0.44\textwidth]{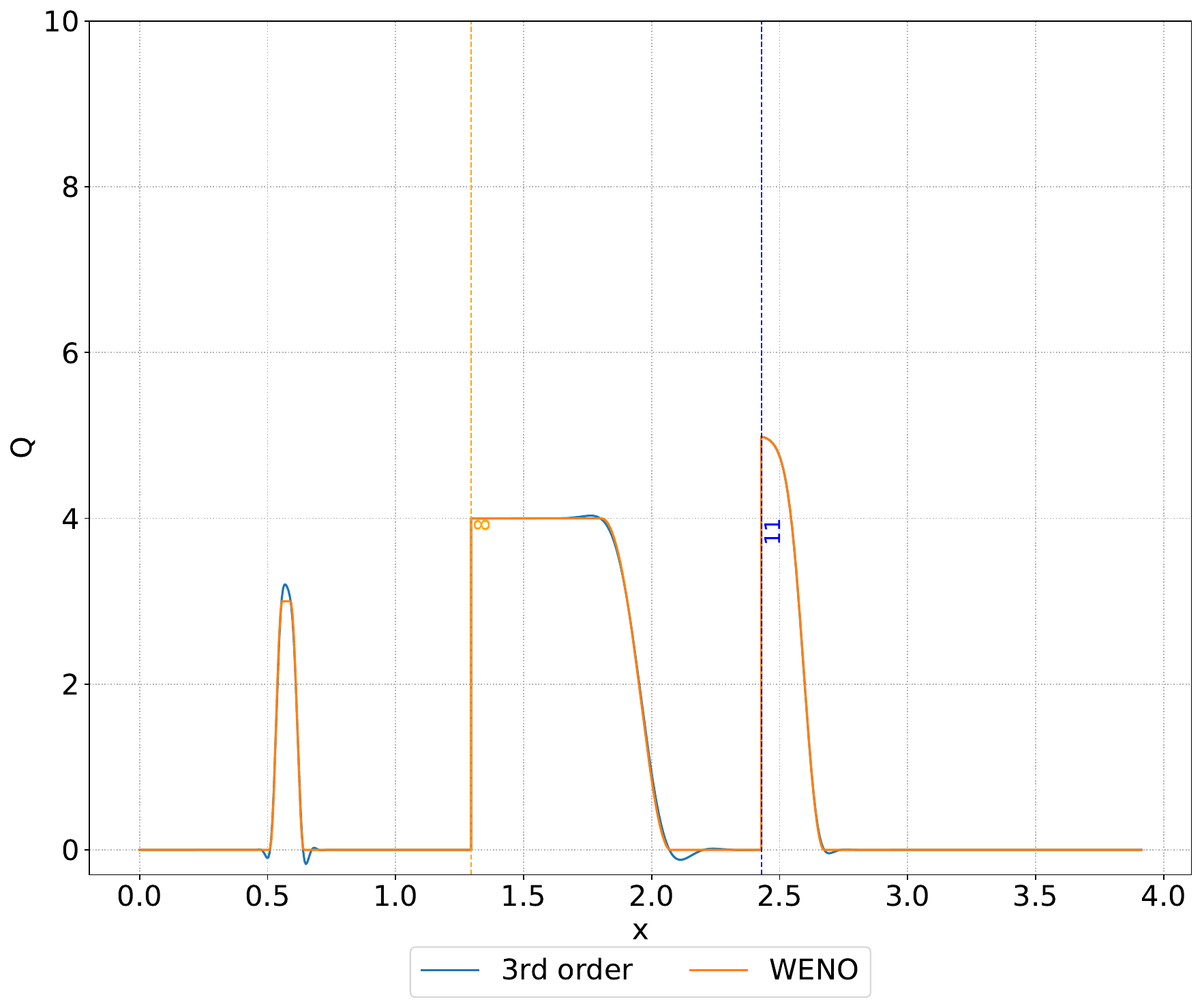}\label{fig:f_bc3_025t_halfCournat}}

  \subfloat[t = 1.0]{\includegraphics[width=0.44\textwidth]{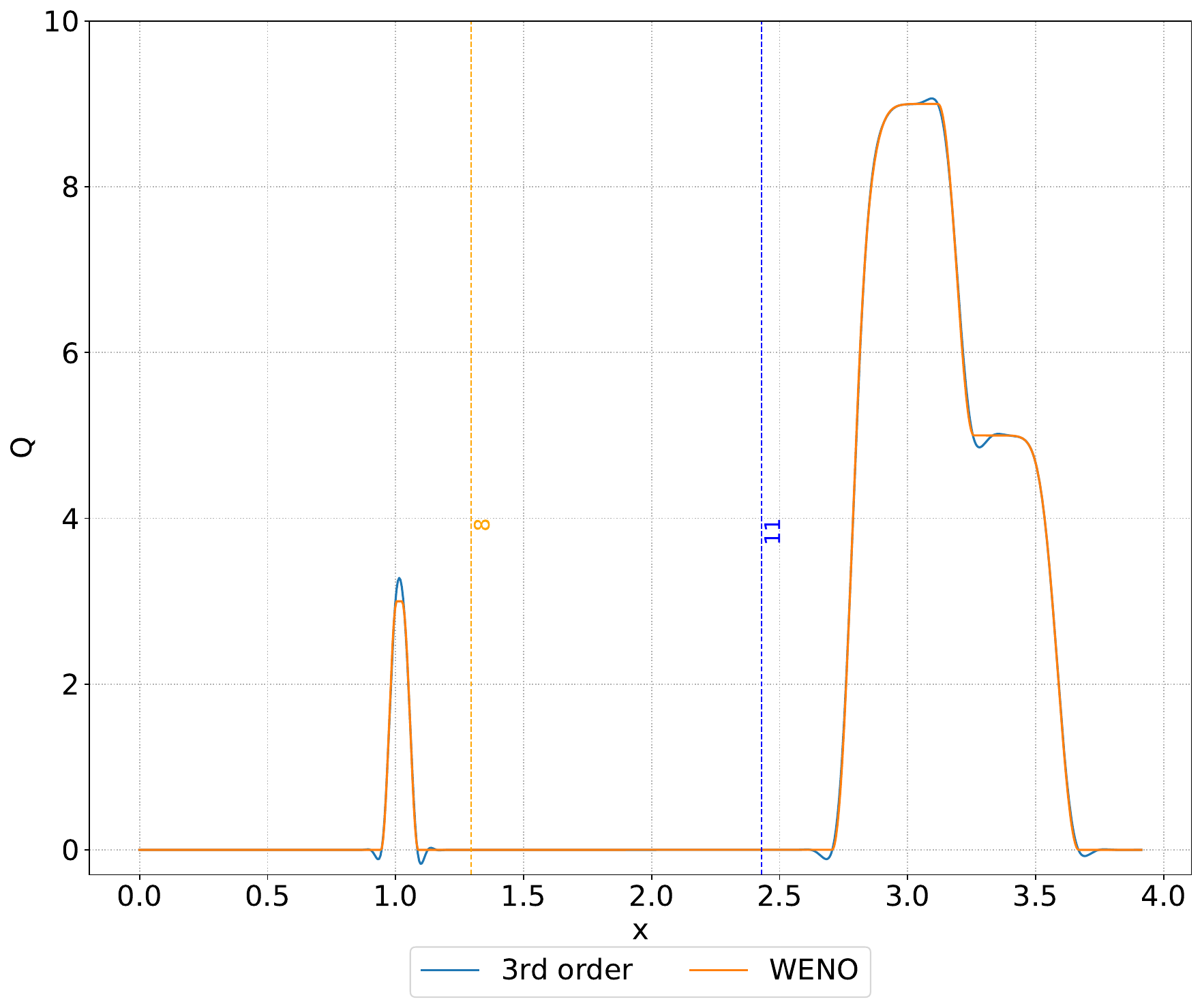}\label{fig:f_bc3_05t_halfCournat}}
  \subfloat[t = 2.0]{\includegraphics[width=0.44\textwidth]{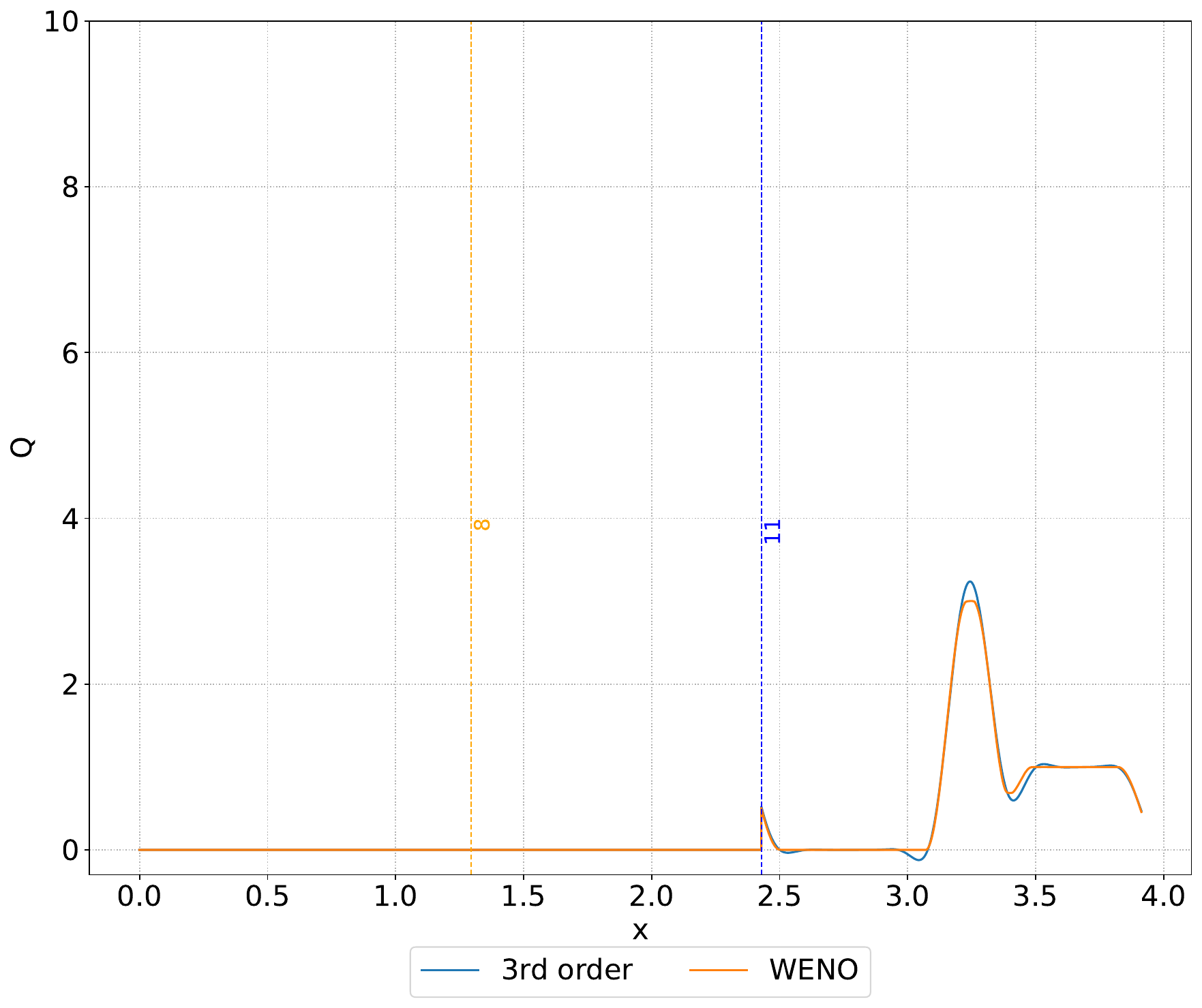}\label{fig:f_bc3_1t_halfCournat}}

  \caption{Flow through the path $BC3$ shown in 4 different moments. Courant numbers on each edge are shown in Table \ref{tab:cournatnumbs_half}.}
  \label{fig:f_way3_timeT_halfCournat}
\end{figure}

\begin{figure}[H]
  \centering
  \subfloat[$N=192$]{\includegraphics[width=0.44\textwidth]{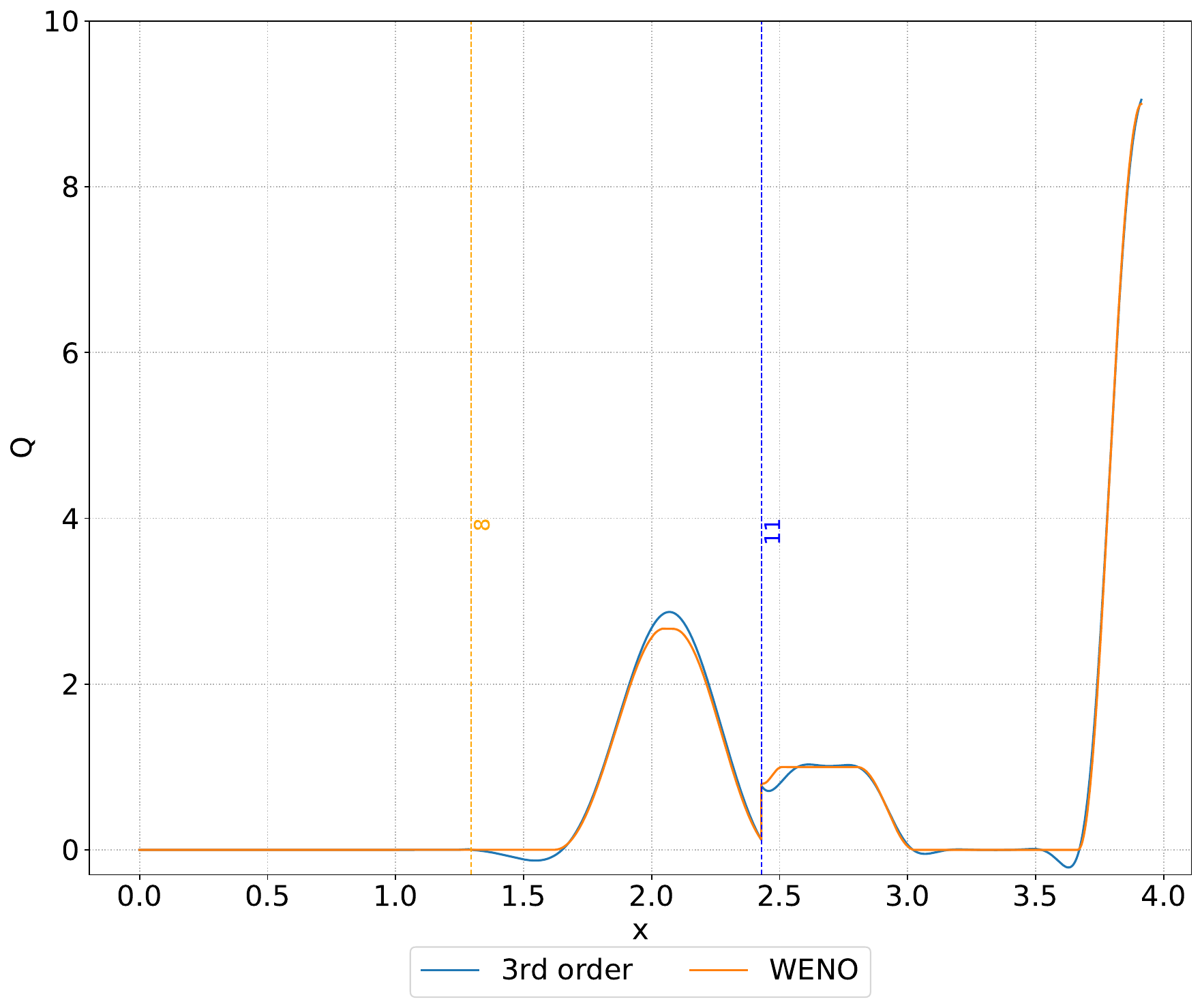}\label{fig:f_T0.7BC1}}
  \subfloat[$N=384$]{\includegraphics[width=0.44\textwidth]{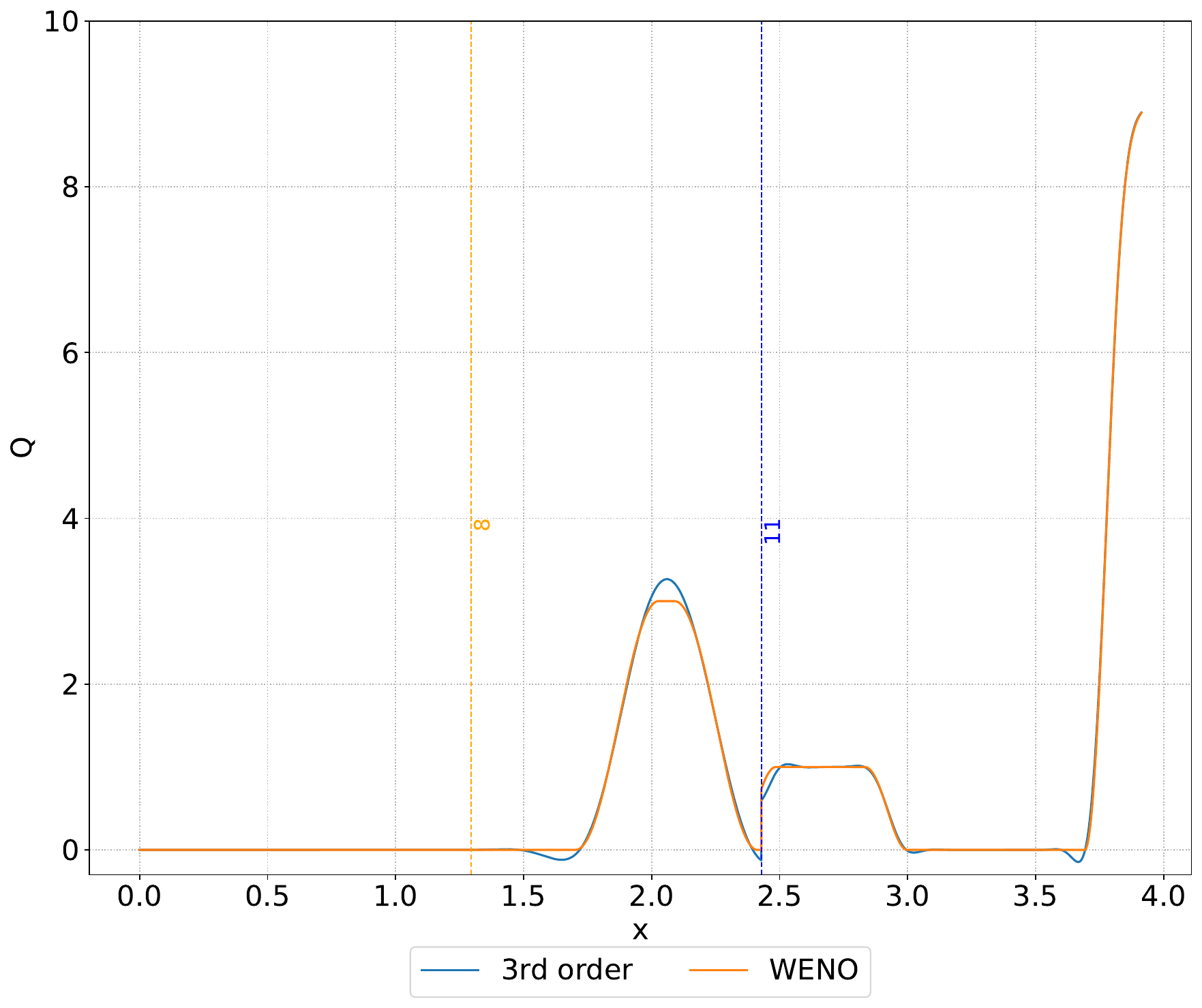}\label{fig:f_T0.7BC3_haflCournt}}
  \caption{Comparison of numerical solutions along the path $BC3$ at time $t = 1.5$ for $N=198$ (the larger Courant numbers) and $N=384$ (the Courant numbers in Table \ref{tab:cournatnumbs_half}. }
  \label{fig:f_all_way_time0.75T_half}
\end{figure}

We comment on the results along the two selected paths. In both cases, clearly, the waves are moving towards the \texttt{OUTFLOW} vertex. The flow profile along both selected paths at different times is shown in the series of Figures \ref{fig:f_way1_timeT_halfCournat} -  \ref{fig:f_all_way_time0.75T_half}. We mark the coupling vertices with dashed lines in these Figures together with the vertex number. Note that we should observe in the solution the waves of the prescribed initial maximal value that are eventually summed up and that can change their lengths due to different speeds on each edge. In Figures, we present the numerical solutions obtained with the 3rd order accurate scheme \eqref{compact-inverse-scheme-fixed-w} with \eqref{parameter-courant-number} and the high-resolution (HR) scheme \eqref{compact-scheme-with-limiters} using single corrector \eqref{corrector}.

Figure \ref{fig:f_way1_timeT_halfCournat} shows the behavior of the solution from the point of view of the path starting at $BC1$. In particular, the wave entering at $BC1$ is visible for all chosen time points, $t=0.25$, $0.5$, $1$, and $2$. For $t\le 1$, it is well approximated by both numerical methods with small over- and undershootings for the third order scheme. Similar conclusions can be drawn for the numerical approximations of the wave entering from $BC2$ that can be observed for $t \in [0.5,2]$.  For $t=1$ and $t=2$ one can observe interactions of the waves. Namely, the wave from $BC1$ interacts with the wave entering from $BC3$ at $t=2$ (see also the description later for the path from $BC3$) where the influence of negative values of the third order scheme for this interaction is clearly visible compared to the HR scheme. Similarly, at $t=1$, the waves from $BC4$ and $BC5$ interact where the influence of unphysical oscillations in the third order scheme is clearly visible.

The difference between the third order scheme and the HR scheme is even more visible in Figure \ref{fig:f_way3_timeT_halfCournat} that shows the transport in the path starting from $BC3$. 
The length of the wave from the boundary condition at $BC3$ can be seen to be shorter in the network due to a slower velocity in this part of the sewer system, so it is harder for numerical methods to approximate it well. The oscillations in the numerical solution obtained with the third order scheme are very clear, with two peaks of the overshooting at $t=0.25$. For $t>0.25$ these two peaks merge, but the magnitude of the overshooting does not seem to decrease even at $t=2$, while the HR scheme preserves the correct maximal value of this wave correctly.

To illustrate the stability of the methods, we present a visual comparison at $t=1.5$ for $N=192$ and identical numbers $I^e$ as in Table \ref{tab:cournatnumbs_half} when Courant numbers take twice as large values. In Figure \ref{fig:f_T0.7BC1} and Figure \ref{fig:f_T0.7BC3_haflCournt} one can observe that the solutions for larger Courant numbers behave qualitatively similarly to those for the smaller ones, with, of course, better precision in the latter case.

For the presented settings, the HR scheme performs very well, effectively suppressing small nonphysical negative values while preserving the flat region of the waves without introducing spurious oscillations. Such behavior seems to be an advantage especially when an interaction between two waves occurs which can be otherwise corrupted by under- or overshooting in approximations of these waves. Note that additional numerical simulations including the time dependent velocity field are presented for this sewer network in \cite{kriskovathesis2025}.

\section{Conclusions}
\label{sec:conclusions}

We have presented an unconditionally bound preserving, high-resolution compact implicit finite difference scheme for solving the advection equation in examples with a predetermined direction of the flow. This type of scheme is suitable for use in network-based models, such as transport in sewage systems, where transport follows the direction from the input vertices to the output vertices. This property allows the scheme to operate as efficiently as if it were explicit, but without any stability restriction on discretization steps. The scheme enables a direct solution by the forward substitution for the linear advection and by solving scalar algebraic equations for each unknown in the case of transport with a nonlinear retardation coefficient. 

The scheme employs the space-time limiter that combines the second-order accurate approximation in space with the third-order approximation in time. It can be derived using the partial inverse Lax-Wendroff procedure, resulting in the finite difference stencil that spans two spatial levels and up to four temporal levels, including the future temporal level. Despite this, it remains efficiently solvable through the nontraditional approach - marching in space instead of time.

To address the nonlinearity in the algebraic equations introduced by the limiter or WENO scheme, we propose the predictor-corrector approach. For chosen examples, the corrector has to be recomputed at most once and only occasionally. We present various benchmark tests and a real-life case study on a sewage network in the city of Revúca. All tests confirmed the expected properties of the scheme. For future development, we plan to consider systems of nonlinear transport equations that can be fully coupled because of the dependence of retardation coefficients on all transported concentrations.

\bibliographystyle{plain}
\bibliography{main}

\begin{thebibliography}{10}

\bibitem{arbogastThirdOrderImplicit2020}
Todd Arbogast, Chieh-Sen Huang, Xikai Zhao, and Danielle~N. King.
\newblock A third order, implicit, finite volume, adaptive {{Runge}}--{{Kutta WENO}} scheme for advection--diffusion equations.
\newblock {\em Comp. Meth. Appl. Mech. Eng.}, 368:113155, 2020.

\bibitem{aroraWellBehavedTVDLimiter1997}
Mohit Arora and Philip~L. Roe.
\newblock A {{Well-Behaved TVD Limiter}} for {{High-Resolution Calculations}} of {{Unsteady Flow}}.
\newblock {\em J. Comp. Phys.}, 132(1):3--11, 1997.

\bibitem{baezaReprintApproximateTaylor2018}
A.~Baeza, S.~Boscarino, P.~Mulet, G.~Russo, and D.~Zor{\'i}o.
\newblock Reprint of: {{Approximate Taylor}} methods for {{ODEs}}.
\newblock {\em Computers \& Fluids}, 169:87--97, 2018.

\bibitem{baezaApproximateImplicitTaylor2020}
Antonio Baeza, Raimund B{\"u}rger, Mar{\'i}a del~Carmen Mart{\'i}, Pep Mulet, and David Zor{\'i}o.
\newblock On approximate implicit {{Taylor}} methods for ordinary differential equations.
\newblock {\em Comp. Appl. Math.}, 39(4):304, 2020.

\bibitem{banik2014swmm5}
{\relax BK}~Banik, C~Di~Cristo, and A~Leopardi.
\newblock {{SWMM5}} toolkit development for pollution source identification in sewer systems.
\newblock {\em Procedia Eng.}, 89:750--757, 2014.

\bibitem{barsukowImplicitActiveFlux2024}
Wasilij Barsukow and Raul Borsche.
\newblock Implicit {{Active Flux Methods}} for {{Linear Advection}}.
\newblock {\em J Sci Comput}, 98(3):52, 2024.

\bibitem{borscheLocalTimeStepping2019}
Raul Borsche, Matthias Eimer, and Norbert Siedow.
\newblock A local time stepping method for thermal energy transport in district heating networks.
\newblock {\em Appl. Math. Comp.}, 353(C):215--229, 2019.

\bibitem{chouchoulisJacobianfreeImplicitMDRK2024}
Jeremy Chouchoulis and Jochen Sch{\"u}tz.
\newblock Jacobian-free implicit {{MDRK}} methods for stiff systems of {{ODEs}}.
\newblock {\em Appl. Num. Math.}, 196:45--61, 2024.

\bibitem{clainHighorderFiniteVolume2011}
S.~Clain, S.~Diot, and R.~Loub{\`e}re.
\newblock A high-order finite volume method for systems of conservation laws---{{Multi-dimensional Optimal Order Detection}} ({{MOOD}}).
\newblock {\em J. Comp. Phys.}, 230(10):4028--4050, 2011.

\bibitem{donatWENOSchemeCharacteristics2024}
R.~Donat, M.~C. Mart{\'i}, and P.~Mulet.
\newblock {{WENO}} scheme on characteristics for the equilibrium dispersive model of chromatography with generalized {{Langmuir}} isotherms.
\newblock {\em Appl. Num. Math.}, 201:247--264, 2024.

\bibitem{duraisamyImplicitSchemeHyperbolic2007}
Karthikeyan Duraisamy and James~D. Baeder.
\newblock Implicit scheme for hyperbolic conservation laws using nonoscillatory reconstruction in space and time.
\newblock {\em SIAM J. Sci. Comput.}, 29(6):2607--2620, 2007.

\bibitem{duraisamyConceptsApplicationTimeLimiters2003}
Karthikeyan Duraisamy, James~D. Baeder, and Jian-Guo Liu.
\newblock Concepts and {{Application}} of {{Time-Limiters}} to {{High Resolution Schemes}}.
\newblock {\em J. Sci. Comp.}, 19(1):139--162, 2003.

\bibitem{eimerImplicitFiniteVolume2022}
Matthias Eimer, Raul Borsche, and Norbert Siedow.
\newblock Implicit finite volume method with a posteriori limiting for transport networks.
\newblock {\em Adv Comput Math}, 48(3):21, 2022.

\bibitem{frolkovicSemianalyticalSolutionsContaminant2006}
Peter Frolkovi{\v c} and Jozef Ka{\v c}ur.
\newblock Semi-analytical solutions of a contaminant transport equation with nonlinear sorption in {{1D}}.
\newblock {\em Comput Geosci}, 10(3):279--290, 2006.

\bibitem{frolkovicSemiimplicitMethodsAdvection2022}
Peter Frolkovi{\v c}, Svetlana Kri{\v s}kov{\'a}, Michaela Rohov{\'a}, and Michal {\v Z}erav{\'y}.
\newblock Semi-implicit methods for advection equations with explicit forms of numerical solution.
\newblock {\em Japan J. Indust. Appl. Math.}, 39(3):843--867, 2022.

\bibitem{frolkovicNumericalSimulationContaminant2016}
Peter Frolkovi{\v c}, Michael Lampe, and Gabriel Wittum.
\newblock Numerical simulation of contaminant transport in groundwater using software tools of r3t.
\newblock {\em Comput. Visual Sci.}, 18(1):17--29, 2016.

\bibitem{frolkovicHighResolutionCompact2023}
Peter Frolkovi{\v c} and Michal {\v Z}erav{\'y}.
\newblock High resolution compact implicit numerical scheme for conservation laws.
\newblock {\em Appl. Math. Comp.}, 442:127720, 2023.

\bibitem{jamesonPositiveSchemesShock1995}
Antony Jameson.
\newblock Positive schemes and shock modelling for compressible flows.
\newblock {\em Int. J. Numer. Methods Fluids}, 20(8-9):743--776, 1995.

\bibitem{kacurSolutionContaminantTransport2005}
Jozef Ka{\v c}ur, Benny Malengier, and Mariana Reme{\v s}{\'i}kov{\'a}.
\newblock Solution of contaminant transport with equilibrium and non-equilibrium adsorption.
\newblock {\em Comp. Meth. Appl. Mech. Eng.}, 194(2):479--489, 2005.

\bibitem{kemmComparativeStudyTVDlimiterswellknown2011}
Friedemann Kemm.
\newblock A comparative study of {{TVD-limiters}}---well-known limiters and an introduction of new ones.
\newblock {\em Int J Numer Meth Fluids}, 67(4):404--440, 2011.

\bibitem{kriskovathesis2025}
Svetlana Kri{\v s}kov{\'a}.
\newblock {\em Numerical Schemes Based on Inverse {{Lax-Wendroff}} Procedure for Advection on Networks}.
\newblock PhD thesis, Slovensk{\'a} technick{\'a} univerzita v Bratislave, Bratislava, 2025.

\bibitem{kuzminDesignGeneralpurposeFlux2006}
D.~Kuzmin.
\newblock On the design of general-purpose flux limiters for finite element schemes. {{I}}. {{Scalar}} convection.
\newblock {\em J. Comp. Phys.}, 219(2):513--531, 2006.

\bibitem{kuzminLocallyBoundpreservingEnriched2020}
Dmitri Kuzmin, Hennes Hajduk, and Andreas Rupp.
\newblock Locally bound-preserving enriched {{Galerkin}} methods for the linear advection equation.
\newblock {\em Computers \& Fluids}, 205:104525, 2020.

\bibitem{lackovaCompactSchemesAdvection2024}
Katar{\'i}na Lackov{\'a} and Peter Frolkovi{\v c}.
\newblock Compact {{Schemes For Advection Equation}}: {{Employing Inverse Lax-wendroff Procedure}}.
\newblock In {\em Algoritmy}, pages 149--158, 2024.

\bibitem{levequeFiniteVolumeMethods2004}
Randall~J. Leveque.
\newblock {\em Finite {{Volume Methods}} for {{Hyperbolic Problems}}}.
\newblock Cambridge UP, 2nd edition, 2004.

\bibitem{maccaSemiimplicitTypeOrderAdaptiveCAT22025}
Emanuele Macca and Sebastiano Boscarino.
\newblock Semi-implicit-{{Type Order-Adaptive CAT2 Schemes}} for {{Systems}} of {{Balance Laws}} with {{Relaxed Source Term}}.
\newblock {\em Commun. Appl. Math. Comput.}, 7(1):151--178, 2025.

\bibitem{maccaAlmostFailsafeAposteriori2024}
Emanuele Macca, Rapha{\"e}l Loub{\`e}re, Carlos Par{\'e}s, and Giovanni Russo.
\newblock An almost fail-safe a-posteriori limited high-order {{CAT}} scheme.
\newblock {\em J. Comp. Phys.}, 498:112650, 2024.

\bibitem{markRiskAnalysesSewer1998}
O.~Mark, C.~Wennberg, T.~{van Kalken}, F.~Rabbi, and B.~Albinsson.
\newblock Risk analyses for sewer systems based on numerical modelling and {{GIS}}.
\newblock {\em Safety Science}, 30(1):99--106, 1998.

\bibitem{michel-dansacTVDMOODSchemesBased2022}
Victor {Michel-Dansac} and Andrea Thomann.
\newblock {{TVD-MOOD}} schemes based on implicit-explicit time integration.
\newblock {\em Appl. Math. Comp.}, 433:127397, 2022.

\bibitem{mohring2021district}
Jan Mohring, Dominik Linn, Matthias Eimer, Markus Rein, and Norbert Siedow.
\newblock District heating networks--dynamic simulation and optimal operation.
\newblock {\em Math. Model. Simul. Optim. Power Eng. Manag.}, pages 303--325, 2021.

\bibitem{pimentel-garciaHighorderIncellDiscontinuous2024}
Ernesto {Pimentel-Garc{\'i}a}, Manuel~J. Castro, Christophe Chalons, and Carlos Par{\'e}s.
\newblock High-order in-cell discontinuous reconstruction path-conservative methods for nonconservative hyperbolic systems--{{DR}}.{{MOOD}} method.
\newblock {\em Num. Meth. Partial Diff. Eq.}, 40(6):e23133, 2024.

\bibitem{puppoQuinpiIntegratingConservation2022}
Gabriella Puppo, Mateo Semplice, and Giuseppe Visconti.
\newblock Quinpi: {{Integrating Conservation Laws}} with {{CWENO Implicit Methods}}.
\newblock {\em Commun. Appl. Math. Comput.}, 2022.

\bibitem{puppoQuinpiIntegratingStiff2024}
Gabriella Puppo, Giuseppe Visconti, and Mateo Semplice.
\newblock Quinpi: {{Integrating Stiff Hyperbolic Systems}} with {{Implicit High Order Finite Volume Schemes}}.
\newblock {\em CiCP}, 36(1):30--70, 2024.

\bibitem{quezadadelunaMaximumPrinciplePreserving2022}
Manuel {Quezada de Luna} and David~I. Ketcheson.
\newblock Maximum {{Principle Preserving Space}} and {{Time Flux Limiting}} for {{Diagonally Implicit Runge}}--{{Kutta Discretizations}} of {{Scalar Convection-diffusion Equations}}.
\newblock {\em J. Sci. Comput.}, 92(3):102, 2022.

\bibitem{roeContributionsNumericalModelling1985}
Philip Roe.
\newblock Some {{Contributions}} to the {{Numerical Modelling}} of {{Discontinuous Flow}}.
\newblock In {\em Large-Scale Computations in Fluid Mechanics}, pages 163--193. American Mathematical Society, 1985.

\bibitem{shuEssentiallyNonoscillatoryWeighted1998}
Chi-Wang Shu.
\newblock Essentially non-oscillatory and weighted essentially non-oscillatory schemes for hyperbolic conservation laws.
\newblock In Bernardo Cockburn, Chi-Wang Shu, Claes Johnson, Eitan Tadmor, and Alfio Quarteroni, editors, {\em Advanced {{Numerical Approximation}} of {{Nonlinear Hyperbolic Equations}}}, Lecture {{Notes}} in {{Mathematics}}, pages 325--432. Springer, Berlin, Heidelberg, 1998.

\bibitem{shuInverseLaxWendroffBoundary2022}
Chi-Wang Shu.
\newblock Inverse {{Lax-Wendroff Boundary Treatment}}: A {{Survey}}.
\newblock {\em CMR}, 38(3):333--350, 2022.

\bibitem{sokacImpactSedimentLayer2021}
Marek Sok{\'a}{\v c} and Yvetta Vel{\'i}skov{\'a}.
\newblock Impact of {{Sediment Layer}} on {{Longitudinal Dispersion}} in {{Sewer Systems}}.
\newblock {\em Water}, 13(22):3168, 2021.

\bibitem{tanInverseLaxWendroffProcedure2010}
Sirui Tan and Chi-Wang Shu.
\newblock Inverse {{Lax-Wendroff}} procedure for numerical boundary conditions of conservation laws.
\newblock {\em J. Comp. Phys.}, 229(21):8144--8166, 2010.

\bibitem{veliskovaInverseTaskPollution2023}
Yvetta Vel{\'i}skov{\'a}, Marek Sok{\'a}{\v c}, and Maryam~Barati Moghaddam.
\newblock Inverse task of pollution spreading -- {{Localization}} of source in extensive open channel network structure.
\newblock {\em J. Hydrol. Hydromech.}, 71(4):475--485, 2023.

\bibitem{zhangReviewTVDSchemes2015}
Di~Zhang, Chunbo Jiang, Dongfang Liang, and Liang Cheng.
\newblock A review on {{TVD}} schemes and a refined flux-limiter for steady-state calculations.
\newblock {\em J. Comp. Phys.}, 302:114--154, 2015.

\bibitem{zorioApproximateLaxWendroff2017}
David Zor{\'i}o, Antonio Baeza, and Pep Mulet.
\newblock An {{Approximate Lax}}--{{Wendroff-Type Procedure}} for {{High Order Accurate Schemes}} for {{Hyperbolic Conservation Laws}}.
\newblock {\em J Sci Comput}, 71(1):246--273, 2017.

\end{thebibliography}

\noindent \section*{In memoriam}
 
\noindent This paper is dedicated to the memory of Prof. Arturo Hidalgo L\'opez
($^*$July 03\textsuperscript{rd} 1966 - $\dagger$August 26\textsuperscript{th} 2024) of the Universidad Politecnica de Madrid, organizer of HONOM 2019 and active participant in many other editions of HONOM.
Our thoughts and wishes go to his wife Lourdes and his sister Mar\'ia Jes\'us, whom he left behind.

\end{document}